\documentclass[final,3p]{elsarticle}
 \usepackage{graphics}
 \usepackage{graphicx}
 \usepackage{epsfig}
\usepackage{amssymb}
 \usepackage{amsthm}
 \usepackage{lineno}
 \usepackage{amsmath}
   \numberwithin{equation}{section}
\usepackage{mathrsfs}

\newtheorem{thm}{Theorem}[section]

\newtheorem{lem}[thm]{Lemma}

\biboptions{sort&compress,square}
\journal{Journal}
\begin{document}
\begin{frontmatter}
\author[rvt1]{Sining Wei}
\ead{weisn835@nenu.edu.cn}
\author[rvt2]{Hongfeng Li}
\ead{lihf728@nenu.edu.cn}
\author[rvt2]{Yong Wang\corref{cor2}}
\ead{wangy581@nenu.edu.cn}

\cortext[cor2]{Corresponding author.}

\address[rvt1]{School of Data Science and Artificial Intelligence, Dongbei University of Finance and Economics, \\
Dalian, 116025, P.R.China}
\address[rvt2]{School of Mathematics and Statistics, Northeast Normal University, Changchun, 130024, P.R.China}

\title{ Conformal perturbations of dirac operators and general Kastler-Kalau-Walze type theorems for even dimensional manifolds with boundary}
\begin{abstract}
In this paper, we establish the proof of  general Kastler-Kalau-Walze type theorems
for conformal perturbations of dirac Operators on even dimensional compact manifolds with (respectively without) boundary.
\end{abstract}
\begin{keyword} Conformal perturbations of dirac Operators; noncommutative residue; Lichnerowicz type formulas; Kastler-Kalau-Walze type theorem.\\

\end{keyword}
\end{frontmatter}
\section{Introduction}
\label{1}
The noncommutative residue plays a significant role in noncommutative geometry, which has been extensively studied by geometers \cite{Gu,Wo}.
Adler discovered the noncommutative residue for one-dimensional manifolds in \cite{MA}, where he explored the geometric aspects of nonlinear partial differential equations.
Wodzicki introduced the noncommutative residue for arbitrary closed compact $n$-dimensional manifolds in [2] using the theory of zeta functions of elliptic pseudodifferential operators.
In \cite{Co1}, Connes used the noncommutative residue to derive a conformal 4-dimensional Polyakov action analogy.
Furthermore, Connes claimed the noncommutative residue of the square of the inverse of the Dirac operator was proportioned to the Einstein-Hilbert action in \cite{Co2}.
Kastler provided a brute-force proof of this theorem in \cite{Ka}, while Kalau and Walze proved it in the normal coordinates system simultaneously in \cite{KW}.  Moreover, Ackermann proved that
the Wodzicki residue of the square of the inverse of the
Dirac operator ${\rm  Wres}(D^{-2})$ in turn is essentially the second coefficient
of the heat kernel expansion of $D^{2}$ in \cite{Ac}.

On the other hand, Fedosov etc. definied a noncommutative residue
on Boutet de Monvel's algebra and proved that it was a unique continuous trace in \cite{FGLS}. Schrohe established a relationship between the Dixmier trace and the noncommutative residue for manifolds with boundary in \cite{S}.
In \cite{Wa3,Wa4}, Wang computed $\widetilde{{\rm Wres}}[\pi^+D^{-1}\circ\pi^+D^{-1}]$ and $\widetilde{{\rm Wres}}[\pi^+D^{-2}\circ\pi^+D^{-2}]$, where the two operators are symmetric, in these cases the boundary term vanished. But for $\widetilde{{\rm Wres}}[\pi^+D^{-1}\circ\pi^+D^{-3}]$, J. Wang and Y. Wang got a nonvanishing boundary term \cite{Wa5}, and give a theoretical explanation for gravitational action on boundary. In others words, Wang provided a kind of method to study the Kastler-Kalau-Walze type theorem for manifolds with boundary.

In \cite{AAA}, Wang established a Kastler-Kalau-Walze type theorem for perturbations of Dirac operators on compact manifolds with (respectively without) boundary.
In \cite{WSW}, Wei and Wang  establish two Kastler-Kalau-Walze type theorems for conformal perturbations of modified
Novikov Operators on 4-dimensional and 6-dimensional compact manifolds with(respectively without) boundary.
In \cite{WJ5}, J. Wang and Y. Wang computed $\widetilde{{\rm Wres}}[(\pi^+D^{-2})\circ(\pi^+D^{-n+2})]$ for manifolds with any dimension and boundary, and established a general Kastler-Kalau-Walze type theorem.
{\bf The motivation of this paper} is to establish the proof of general Kastler-Kalau-Walze type theorems for conformal perturbations of dirac Operators on even dimensional compact manifolds with (respectively without) boundary.
In this paper, the leading symbol of dirac operator is $ic(\xi)$. At the moment,
the leading symbol of conformal perturbations
of dirac operators is not $ic(\xi)$, which motivates the study of the residue of conformal perturbations
of dirac operators.
That is, we want to compute ${\rm Wres}[\pi^+P_1\circ\pi^+P_2]$, where orders of $P_1,P_2$ are $a_1,a_2$ and $-a_1-a_2+2=m$ for even dimensional manifolds with boundary. Motivated by \cite{WJ5,WSW}, we compute the generalized noncommutative residue $\widetilde{{\rm Wres}}\bigg[\pi^+(fD^{-1}f^{-1}D^{-1})\circ\pi^+\Big((fD^{-1}f^{-1}D^{-1})^{n}\Big)\bigg]$ and $\widetilde{{\rm Wres}}\bigg[\pi^+(fD^{-1})\circ\pi^+\Big((f^{-1}D^{-1})\cdot(fD^{-1}\cdot f^{-1}D^{-1})^{n}\Big)\bigg]$ on even dimensional manifolds. Our main theorems are as follows.

\begin{thm}\label{thm1.1}
Let $M$ be an $m=2n+4$ dimensional oriented
compact spin manifold with boundary $\partial M$, then we get the following equality:
\begin{eqnarray}
\label{1111}
&&\widetilde{{\rm Wres}}\bigg[\pi^+(fD^{-1}f^{-1}D^{-1})\circ\pi^+\Big((fD^{-1}f^{-1}D^{-1})^{n}\Big)\bigg]\nonumber\\
&=&\frac{(2\pi)^{2n+6}}{(2n+4)!}\int_{M}2^{2n+6}
\bigg\{-\frac{1}{12}s-2f^{-1}\Delta(f)-f^{-2}\Big[|\mathrm{grad}_{M}f|^2+2\Delta(f)\Big]\bigg\}d{\rm Vol_{M}}\nonumber\\
&&+\int_{\partial M}\bigg\{\frac{2^{n}h'(0)\pi ni}{(n+3)!}{\rm Vol}(S_{2n+2})
Q_0\bigg\}d{\rm Vol}_{\partial_{M}},
\end{eqnarray}
where $Q_0$ are defined in (\ref{fum9.1}).
\end{thm}
\begin{thm}\label{thm1.2}
Let $M$ be an $m=2n+4$ dimensional oriented
compact spin manifold with boundary $\partial M$, then we get the following equality:
\begin{eqnarray}
\label{1b2773}
&&\widetilde{{\rm Wres}}\bigg[\pi^+(fD^{-1})\circ\pi^+\Big((f^{-1}D^{-1})\cdot(fD^{-1}\cdot f^{-1}D^{-1})^{n}\Big)\bigg]\nonumber\\
&=& \frac{(2\pi)^{2n+6}}{(2n+4)!}\int_{M}2^{2n+6}
\bigg\{-\frac{1}{12}s-2f^{-1}\Delta(f)-f^{-2}\Big[|\mathrm{grad}_{M}f|^2+2\Delta(f)\Big]\bigg\}d{\rm Vol_{M}}\nonumber\\
&&+\int_{\partial M}\Bigg\{\frac{(-1)^{n}h'(0)\pi }{3\times2^{n+6}(3+n)!}Y_0
+\frac{(-1)^{n}(n+1)f^{-1}\partial_{x_{n}}(f)\pi}{2^{n+2}}Y_1
+\bigg[\frac{(i-1)f\cdot\partial_{x_{n}}(f^{-1})}{2^{n+3}}\nonumber\\
&&-\frac{\partial_{x_{n}}(f)(1+i)}{2^{n+2}} \bigg]Y_2\Bigg\}d{\rm Vol}_{\partial_{M}},
\end{eqnarray}
where $Y_0,~Y_1,~Y_2$ are defined in (\ref{fum11.1}).
\end{thm}

The paper is organized in the following way. In Section \ref{section:2}, we review
some basic formulas related to Boutet de
Monvel's calculus and the definition of the noncommutative residue for manifolds with boundary. In Section \ref{section:3}, we prove the general Kastler-Kalau-Walze type theorem for $\widetilde{{\rm Wres}}\bigg[\pi^+(fD^{-1}f^{-1}D^{-1})\circ\pi^+\Big((fD^{-1}f^{-1}D^{-1})^{n}\Big)\bigg]$ on even dimensional manifolds with boundary. In Section \ref{section:4},  we prove the general Kastler-Kalau-Walze type theorem $\widetilde{{\rm Wres}}\bigg[\pi^+(fD^{-1})\circ\pi^+\Big((f^{-1}D^{-1})\cdot(fD^{-1}\cdot f^{-1}D^{-1})^{n}\Big)\bigg]$ on even dimensional manifolds with boundary.

\section{Boutet de
Monvel's calculus and the definition of the noncommutative residue}
\label{section:2}
 In this section, we recall some basic facts and formulas about Boutet de
Monvel's calculus and the definition of the noncommutative residue for manifolds with boundary which will be used in the following. For more details, see Section 2 in \cite{Wa3}.\\
 \indent Let $M$ be a 4-dimensional compact oriented manifold with boundary $\partial M$.
We assume that the metric $g^{TM}$ on $M$ has the following form near the boundary,
\begin{equation}
\label{b1}
g^{M}=\frac{1}{h(x_{n})}g^{\partial M}+dx _{n}^{2},
\end{equation}
where $g^{\partial M}$ is the metric on $\partial M$ and $h(x_n)\in C^{\infty}([0, 1)):=\{\widehat{h}|_{[0,1)}|\widehat{h}\in C^{\infty}((-\varepsilon,1))\}$ for
some $\varepsilon>0$ and $h(x_n)$ satisfies $h(x_n)>0$, $h(0)=1$ where $x_n$ denotes the normal directional coordinate. Let $U\subset M$ be a collar neighborhood of $\partial M$ which is diffeomorphic with $\partial M\times [0,1)$. By the definition of $h(x_n)\in C^{\infty}([0,1))$
and $h(x_n)>0$, there exists $\widehat{h}\in C^{\infty}((-\varepsilon,1))$ such that $\widehat{h}|_{[0,1)}=h$ and $\widehat{h}>0$ for some
sufficiently small $\varepsilon>0$. Then there exists a metric $g'$ on $\widetilde{M}=M\bigcup_{\partial M}\partial M\times
(-\varepsilon,0]$ which has the form on $U\bigcup_{\partial M}\partial M\times (-\varepsilon,0 ]$
\begin{equation}
\label{b2}
g'=\frac{1}{\widehat{h}(x_{n})}g^{\partial M}+dx _{n}^{2} ,
\end{equation}
such that $g'|_{M}=g$. We fix a metric $g'$ on the $\widetilde{M}$ such that $g'|_{M}=g$.

Let the Fourier transformation $F'$  be
\begin{equation*}
F':L^2({\bf R}_t)\rightarrow L^2({\bf R}_v);~F'(u)(v)=\int_\mathbb{R} e^{-ivt}u(t)dt
\end{equation*}
and let
\begin{equation*}
r^{+}:C^\infty ({\bf R})\rightarrow C^\infty (\widetilde{{\bf R}^+});~ f\rightarrow f|\widetilde{{\bf R}^+};~
\widetilde{{\bf R}^+}=\{x\geq0;x\in {\bf R}\}.
\end{equation*}
\indent We define $H^+=F'(\Phi(\widetilde{{\bf R}^+}));~ H^-_0=F'(\Phi(\widetilde{{\bf R}^-}))$ which satisfies
$H^+\bot H^-_0$, where $\Phi(\widetilde{{\bf R}^+}) =r^+\Phi({\bf R})$, $\Phi(\widetilde{{\bf R}^-}) =r^-\Phi({\bf R})$ and $\Phi({\bf R})$
denotes the Schwartz space. We have the following
 property: $h\in H^+~$ (respectively $H^-_0$) if and only if $h\in C^\infty({\bf R})$ which has an analytic extension to the lower (respectively upper) complex
half-plane $\{{\rm Im}\xi<0\}$ (respectively $\{{\rm Im}\xi>0\})$ such that for all nonnegative integer $l$,
 \begin{equation*}
\frac{d^{l}h}{d\xi^l}(\xi)\sim\sum^{\infty}_{k=1}\frac{d^l}{d\xi^l}(\frac{c_k}{\xi^k}),
\end{equation*}
as $|\xi|\rightarrow +\infty,{\rm Im}\xi\leq0$ (respectively ${\rm Im}\xi\geq0)$ and where $c_k\in\mathbb{C}$ are some constants.\\
 \indent Let $H'$ be the space of all polynomials and $H^-=H^-_0\bigoplus H';~H=H^+\bigoplus H^-.$ Denote by $\pi^+$ (respectively $\pi^-$) the
 projection on $H^+$ (respectively $H^-$). Let $\widetilde H=\{$rational functions having no poles on the real axis$\}$. Then on $\tilde{H}$,
 \begin{equation}
 \label{b3}
\pi^+h(\xi_0)=\frac{1}{2\pi i}\lim_{u\rightarrow 0^{-}}\int_{\Gamma^+}\frac{h(\xi)}{\xi_0+iu-\xi}d\xi,
\end{equation}
where $\Gamma^+$ is a Jordan closed curve
included ${\rm Im}(\xi)>0$ surrounding all the singularities of $h$ in the upper half-plane and
$\xi_0\in {\bf R}$. In our computations, we only compute $\pi^+h$ for $h$ in $\widetilde{H}$. Similarly, define $\pi'$ on $\tilde{H}$,
\begin{equation}
\label{b4}
\pi'h=\frac{1}{2\pi}\int_{\Gamma^+}h(\xi)d\xi.
\end{equation}
So $\pi'(H^-)=0$. For $h\in H\bigcap L^1({\bf R})$, $\pi'h=\frac{1}{2\pi}\int_{{\bf R}}h(v)dv$ and for $h\in H^+\bigcap L^1({\bf R})$, $\pi'h=0$.\\
\indent An operator of order $m\in {\bf Z}$ and type $d$ is a matrix\\
$$\widetilde{A}=\left(\begin{array}{lcr}
  \pi^+P+G  & K  \\
   T  &  \widetilde{S}
\end{array}\right):
\begin{array}{cc}
\   C^{\infty}(M,E_1)\\
 \   \bigoplus\\
 \   C^{\infty}(\partial{M},F_1)
\end{array}
\longrightarrow
\begin{array}{cc}
\   C^{\infty}(M,E_2)\\
\   \bigoplus\\
 \   C^{\infty}(\partial{M},F_2)
\end{array},
$$
where $M$ is a manifold with boundary $\partial M$ and
$E_1,E_2$~ (respectively $F_1,F_2$) are vector bundles over $M~$ (respectively $\partial M
$).~Here,~$P:C^{\infty}_0(\Omega,\overline {E_1})\rightarrow
C^{\infty}(\Omega,\overline {E_2})$ is a classical
pseudodifferential operator of order $m$ on $\Omega$, where
$\Omega$ is a collar neighborhood of $M$ and
$\overline{E_i}|M=E_i~(i=1,2)$. $P$ has an extension:
$~{\cal{E'}}(\Omega,\overline {E_1})\rightarrow
{\cal{D'}}(\Omega,\overline {E_2})$, where
${\cal{E'}}(\Omega,\overline {E_1})~({\cal{D'}}(\Omega,\overline
{E_2}))$ is the dual space of $C^{\infty}(\Omega,\overline
{E_1})~(C^{\infty}_0(\Omega,\overline {E_2}))$. Let
$e^+:C^{\infty}(M,{E_1})\rightarrow{\cal{E'}}(\Omega,\overline
{E_1})$ denote extension by zero from $M$ to $\Omega$ and
$r^+:{\cal{D'}}(\Omega,\overline{E_2})\rightarrow
{\cal{D'}}(\Omega, {E_2})$ denote the restriction from $\Omega$ to
$X$, then define
$$\pi^+P=r^+Pe^+:C^{\infty}(M,{E_1})\rightarrow {\cal{D'}}(\Omega,
{E_2}).$$ In addition, $P$ is supposed to have the
transmission property; this means that, for all $j,k,\alpha$, the
homogeneous component $p_j$ of order $j$ in the asymptotic
expansion of the
symbol $p$ of $P$ in local coordinates near the boundary satisfies:\\
$$\partial^k_{x_n}\partial^\alpha_{\xi'}p_j(x',0,0,+1)=
(-1)^{j-|\alpha|}\partial^k_{x_n}\partial^\alpha_{\xi'}p_j(x',0,0,-1),$$
then $\pi^+P:C^{\infty}(M,{E_1})\rightarrow C^{\infty}(M,{E_2})$. Let $G$, $T$ be respectively the singular Green operator
and the trace operator of order $m$ and type $d$. Let $K$ be a
potential operator and $S$ be a classical pseudodifferential
operator of order $m$ along the boundary. Denote by $B^{m,d}$ the collection of all operators of
order $m$
and type $d$,  and $\mathcal{B}$ is the union over all $m$ and $d$.\\
\indent Recall that $B^{m,d}$ is a Fr\'{e}chet space. The composition
of the above operator matrices yields a continuous map:
$B^{m,d}\times B^{m',d'}\rightarrow B^{m+m',{\rm max}\{
m'+d,d'\}}.$ Write $$\widetilde{A}=\left(\begin{array}{lcr}
 \pi^+P+G  & K \\
 T  &  \widetilde{S}
\end{array}\right)
\in B^{m,d},
 \widetilde{A}'=\left(\begin{array}{lcr}
\pi^+P'+G'  & K'  \\
 T'  &  \widetilde{S}'
\end{array} \right)
\in B^{m',d'}.$$\\
The composition $\widetilde{A}\widetilde{A}'$ is obtained by multiplication of the matrices (For more details see \cite{Ka}).
For
example $\pi^+P\circ G'$ and $G\circ G'$ are singular Green
operators of type $d'$ and
$$\pi^+P\circ\pi^+P'=\pi^+(PP')+L(P,P').$$
Here $PP'$ is the usual
composition of pseudodifferential operators and $L(P,P')$ called
leftover term is a singular Green operator of type $m'+d$. For our case, $P,~P'$ are classical pseudo differential operators, in other words $\pi^+P\in \mathcal{B}^{\infty}$ and $\pi^+P'\in \mathcal{B}^{\infty}$ .\\
\indent Let $M$ be a $n$-dimensional compact oriented manifold with boundary $\partial M$.
Denote by $\mathcal{B}$ the Boutet de Monvel's algebra. We recall that the main theorem in \cite{FGLS,Wa3}.
\begin{thm}\label{thm2.1}{\rm\cite{FGLS}}{\bf(Fedosov-Golse-Leichtnam-Schrohe)}
 Let $M$ and $\partial M$ be connected, ${\rm dim}M=n\geq3$, and let $\widetilde{S}$ (respectively $\widetilde{S}'$) be the unit sphere about $\xi$ (respectively $\xi'$) and $\sigma(\xi)$ (respectively $\sigma(\xi')$) be the corresponding canonical
$n-1$ (respectively $(n-2)$) volume form.
 Set $\widetilde{A}=\left(\begin{array}{lcr}\pi^+P+G &   K \\
T &  \widetilde{S}    \end{array}\right)$ $\in \mathcal{B}$ , and denote by $p$, $b$ and $s$ the local symbols of $P,G$ and $\widetilde{S}$ respectively.
 Define:
 \begin{align}
{\rm{\widetilde{Wres}}}(\widetilde{A})&=\int_X\int_{\bf \widetilde{ S}}{\rm{trace}}_E\left[p_{-n}(x,\xi)\right]\sigma(\xi)dx \nonumber\\
&+2\pi\int_ {\partial X}\int_{\bf \widetilde{S}'}\left\{{\rm trace}_E\left[({\rm{tr}}b_{-n})(x',\xi')\right]+{\rm{trace}}
_F\left[s_{1-n}(x',\xi')\right]\right\}\sigma(\xi')dx',
\end{align}
where ${\rm{\widetilde{Wres}}}$ denotes the noncommutative residue of an operator in the Boutet de Monvel's algebra.\\
Then~~ a) ${\rm \widetilde{Wres}}([\widetilde{A},B])=0 $, for any
$\widetilde{A},B\in\mathcal{B}$;~~ b) It is the unique continuous trace on
$\mathcal{B}/\mathcal{B}^{-\infty}$.
\end{thm}

\begin{thm}\label{thm2.3}
{\rm\cite{WJ6}} For even m-dimensional compact spin manifolds without boundary, the following equality holds:
 \begin{eqnarray}
{\rm Wres}\bigg[fDf^{-1}D^{-1} \bigg]^{\frac{m-2}{2}}
=\frac{(2\pi)^{\frac{m}{2}}}{(\frac{m}{2}-2)!}\int_{M}{\rm{trace}}
\bigg\{-\frac{1}{12}s-2f^{-1}\Delta(f)-f^{-2}[|\mathrm{grad}_{M}f|^2+2\Delta(f)]\bigg \}dvol_{M},
\end{eqnarray}
where s is the scaler curvature.
\end{thm}

\section{ The noncommutative residue $\widetilde{{\rm Wres}}\bigg[\pi^+(fD^{-1}f^{-1}D^{-1})\circ\pi^+\Big((fD^{-1}f^{-1}D^{-1})^{n}\Big)\bigg]$ on even dimensional manifolds with boundary}
\label{section:3}
Firstly, we recall the definition of the Dirac operator. Let $M$ be an $m=2n+4$ dimensional oriented compact spin Riemannian manifold with a Riemannian metric $g^{M}$ and let $\nabla^L$ be the Levi-Civita connection about $g^{M}$.

Set  $\widetilde{e}_{m}=\frac{\partial}{\partial x_{m}}$, $\widetilde{e}_{j}=\sqrt{h(x_{m})}e_{j}~~(1\leq j \leq m-1)$, where  $\{e_{1},\cdots,e_{m-1}\}$ are orthonormal basis of $T\partial_{M}$.
In the local coordinates $\{x_i; 1\leq i\leq m\}$ and the fixed orthonormal frame $\{\widetilde{e}_1,\cdots,\widetilde{e}_m\}$, the connection matrix $(\omega_{s,t})$ is defined by
\begin{equation}
\label{a2}
\nabla^L(\widetilde{e}_1,\cdots,\widetilde{e}_m)= (\widetilde{e}_1,\cdots,\widetilde{e}_m)(\omega_{s,t}).
\end{equation}
\indent Let $c(\widetilde{e}_i)$ denotes the Clifford action, which satisfies
\begin{align}
\label{a4}
&c(\widetilde{e}_i)c(\widetilde{e}_j)+c(\widetilde{e}_j)c(\widetilde{e}_i)=-2g^{M}(\widetilde{e}_i,\widetilde{e}_j).
\end{align}
In \cite{Y}, the Dirac operator is given
\begin{align}
\label{a5}
D=\sum^m_{i=1}c(\widetilde{e}_i)\bigg[\widetilde{e}_i-\frac{1}{4}\sum_{s,t}\omega_{s,t}
(\widetilde{e}_i)
c(\widetilde{e}_s)c(\widetilde{e}_t)\bigg].
\end{align}
Set a Clifford action $c(X)$ on $M$ and $X=\sum\limits_{\alpha=1}^ma_{\alpha}\widetilde{e}_\alpha=X^T+X_m\partial_{x_m}=\sum\limits_{j=1}^mX_j\partial_j$ is a vector field. We define $\nabla_X^{S(TM)}:=X+\frac{1}{4}\sum\limits_{ij}\langle\nabla_X^L{\widetilde{e}_i},\widetilde{e}_j\rangle c(\widetilde{e}_i)c(\widetilde{e}_j)$, which is a spin connection, where $L(X)=\frac{1}{4}\sum\limits_{ij}\langle\nabla_X^L{\widetilde{e}_i},\widetilde{e}_j\rangle c(\widetilde{e}_i)c(\widetilde{e}_j)$. And let $g^{ij}=g(dx_{i},dx_{j})$, $\xi=\sum\limits_{k}\xi_{j}dx_{j}$ and $\nabla^L_{\partial_{i}}\partial_{j}=\sum\limits_{k}\Gamma_{ij}^{k}\partial_{k}$,  we denote that
\begin{align}
&\sigma_{i}=-\frac{1}{4}\sum_{s,t}\omega_{s,t}
(\widetilde{e}_i)c(\widetilde{e}_i)c(\widetilde{e}_s)c(\widetilde{e}_t)
;~~~\xi^{j}=g^{ij}\xi_{i};~~~~\Gamma^{k}=g^{ij}\Gamma_{ij}^{k};~~~~\sigma^{j}=g^{ij}\sigma_{i}.
\end{align}
Then by \cite{Wa3} and $\sigma(\partial_{x_j})=i\xi_j$, we have the following lemmas.
\begin{lem}\label{lem2} The following identities hold:
\begin{align}
&\sigma_1(D)=ic(\xi); \nonumber\\
&\sigma_0(D)=-\frac{1}{4}\sum_{i,s,t}\omega_{s,t}(\widetilde{e}_i)c(\widetilde{e}_i)c(\widetilde{e}_s)c(\widetilde{e}_t)\nonumber\\
&\sigma_{0}(\nabla_X^{S(TM)})=L(X);\nonumber\\
&\sigma_{1}(\nabla_X^{S(TM)})=i\sum_{j=1}^nX_j\xi_j.\nonumber
\end{align}
\end{lem}

By the composition formula of pseudodifferential operators, we have
\begin{lem}\label{lem3} The following identities hold:
\begin{align}
\label{ab22}
&\sigma_{-1}({D}^{-1})=\frac{\sqrt{-1}c(\xi)}{|\xi|^2};\nonumber\\
&\sigma_{-2}(fD^{-1}f^{-1}D^{-1})=\sigma_{-2}(D^{-2})=|\xi|^{-2};\nonumber\\
&\sigma_{-2n}\Big[(fD^{-1}f^{-1}D^{-1})^{n}\Big]=\sigma_{-2n}(D^{-2n})=|\xi|^{-2n};\nonumber\\
&\sigma_{-2n-1}(D^{-2n-1})=\sqrt{-1}c(\xi)|\xi|^{-2n-2};\nonumber\\
&\sigma_{-2}({D}^{-1})=\frac{c(\xi)\sigma_{0}(D)c(\xi)}{|\xi|^4}+\frac{c(\xi)}{|\xi|^6}\sum_jc(dx_j)
\Big[\partial_{x_j}(c(\xi))|\xi|^2-c(\xi)\partial_{x_j}(|\xi|^2)\Big]\nonumber\\
&\sigma_{-3}(D^{-2})=-\sqrt{-1}|\xi|^{-4}\xi_k(\Gamma^k-2\sigma^k)-\sqrt{-1}|\xi|^{-6}2\xi^j\xi_\alpha\xi_\beta\partial_jg^{\alpha\beta};\nonumber\\
&\sigma_{-2n-1}\bigg((fD^{-1}f^{-1}D^{-1})^{n}\bigg)=n\cdot\sigma_{2}^{(1-n)}\sigma_{-3}(fD^{-1}f^{-1}D^{-1})-
          i\cdot \sum_{k=0}^{n-2}\partial_{\xi_{\mu}}
\sigma_{2}^{(1-n+k)}\partial_{x_{\mu}}\sigma_{2}^{-1}\big(\sigma_{2}^{-1}\big)^{k},
\end{align}
where $\sigma_{2}=(1+\xi_{m}^2)^2.$
\end{lem}
Since $\Theta$ is a global form on $\partial M$, so for any fixed point $x_{0}\in\partial M$, we can choose the normal coordinates
$U$ of $x_{0}$ in $\partial M$(not in $M$) and compute $\Theta(x_{0})$ in the coordinates $\widetilde{U}=U\times [0,1)$ and the metric
$\frac{1}{h(x_{m})}g^{\partial M}+dx _{m}^{2}$. The dual metric of $g^{M}$ on $\widetilde{U}$ is
$h(x_{m})g^{\partial M}+dx _{m}^{2}.$ Write
$g_{ij}^{M}=g^{M}(\frac{\partial}{\partial x_{i}},\frac{\partial}{\partial x_{j}})$;
$g^{ij}_{M}=g^{M}(d x_{i},dx_{j})$, then

\begin{equation*}
[g_{i,j}^{M}]=
\begin{bmatrix}\frac{1}{h( x_{m})}[g_{i,j}^{\partial M}]&0\\0&1\end{bmatrix};\quad
[g^{i,j}_{M}]=\begin{bmatrix} h( x_{m})[g^{i,j}_{\partial M}]&0\\0&1\end{bmatrix},
\end{equation*}
and
\begin{equation*}
\partial_{x_{s}} g_{ij}^{\partial M}(x_{0})=0,\quad 1\leq i,j\leq m-1;\quad g_{i,j}^{M}(x_{0})=\delta_{ij}.
\end{equation*}

Let $\{e_{1},\cdots, e_{m-1}\}$ be an orthonormal frame field in $U$ about $g^{\partial M}$ which is parallel along geodesics and
$e_{i}=\frac{\partial}{\partial x_{i}}(x_{0})$, then $\{\widetilde{e_{1}}=\sqrt{h(x_{m})}e_{1}, \cdots,
\widetilde{e_{m-1}}=\sqrt{h(x_{n})}e_{m-1},\widetilde{e_{m}}=dx_{m}\}$ is the orthonormal frame field in $\widetilde{U}$ about $g^{M}.$
Locally $S(TM)|\widetilde{U}\cong \widetilde{U}\times\wedge^{*}_{C}(\frac{m}{2}).$ Let $\{f_{1},\cdots,f_{m}\}$ be the orthonormal basis of
$\wedge^{*}_{C}(\frac{m}{2})$. Take a spin frame field $\sigma: \widetilde{U}\rightarrow Spin(M)$ such that
$\pi\sigma=\{\widetilde{e_{1}},\cdots, \widetilde{e_{m}}\}$ where $\pi: Spin(M)\rightarrow O(M)$ is a double covering, then
$\{[\sigma, f_{i}], 1\leq i\leq m\}$ is an orthonormal frame of $S(TM)|_{\widetilde{U}}$. In the following, since the global form $\Theta$
is independent of the choice of the local frame, so we can compute $\texttt{tr}_{S(TM)}$ in the frame $\{[\sigma, f_{i}], 1\leq i\leq m\}$.
Let $\{\hat{e}_{1},\cdots,\hat{e}_{m}\}$ be the canonical basis of $\mathbb{R}^{m}$ and
$c(\hat{e}_{i})\in cl_{C}(m)\cong Hom(\wedge^{*}_{C}(\frac{m}{2}),\wedge^{*}_{C}(\frac{m}{2}))$ be the Clifford action. Then
\begin{equation*}
c(\widetilde{e_{i}})=[(\sigma,c(\hat{e}_{i}))]; \quad c(\widetilde{e_{i}})[(\sigma, f_{i})]=[\sigma,(c(\hat{e}_{i}))f_{i}]; \quad
\frac{\partial}{\partial x_{i}}=[(\sigma,\frac{\partial}{\partial x_{i}})],
\end{equation*}
then we have $\frac{\partial}{\partial x_{i}}c(\widetilde{e_{i}})=0$ in the above frame. By Lemma 2.2 in \cite{Wa3}, we have

\begin{lem}\label{le:32}
With the metric $g^{M}$ on $M$ near the boundary
\begin{eqnarray}\label{fum100}
\partial_{x_j}(|\xi|_{g^M}^2)(x_0)&=&\left\{
       \begin{array}{c}
        0,  ~~~~~~~~~~ ~~~~~~~~~~ ~~~~~~~~~~~~~{\rm if }~j<m; \\[2pt]
       h'(0)|\xi'|^{2}_{g^{\partial M}},~~~~~~~~~~~~~~~~~~~~~{\rm if }~j=m.
       \end{array}
    \right. \\
\partial_{x_j}[c(\xi)](x_0)&=&\left\{
       \begin{array}{c}
      0,  ~~~~~~~~~~ ~~~~~~~~~~ ~~~~~~~~~~~~~{\rm if }~j<m;\\[2pt]
\partial x_{n}(c(\xi'))(x_{0}), ~~~~~~~~~~~~~~~~~{\rm if }~j=m,
       \end{array}
    \right.
\end{eqnarray}
where $\xi=\xi'+\xi_{m}dx_{m}$.
\end{lem}

In the following, we will compute the residue $\widetilde{{\rm Wres}}\bigg[\pi^+(fD^{-1}f^{-1}D^{-1})\circ\pi^+\Big((fD^{-1}f^{-1}D^{-1})^{n}\Big)\bigg]$ for nonzero
smooth functions $f,~f^{-1}$ on even dimensional oriented
compact spin manifolds with boundary and get a general Kastler-Kalau-Walze
type theorem in this case. By Theorem \ref{thm2.1}, we have
\begin{eqnarray}\label{fum3.7}
&&\widetilde{{\rm Wres}}\bigg[\pi^+(fD^{-1}f^{-1}D^{-1})\circ\pi^+\Big((fD^{-1}f^{-1}D^{-1})^{n}\Big)\bigg]\nonumber\\
&=&\int_{M}\int_{|\xi|=1}{{\rm trace}}_{S(TM)}\big[\sigma_{-n}\big( (fD^{-1}f^{-1}D^{-1})^{n+1}\big)\big]\sigma(\xi)dx+\int_{\partial M}\Phi,
\end{eqnarray}
where
 \begin{eqnarray}\label{fum3.81}
\Phi&=&\int_{|\xi'|=1}\int_{-\infty}^{+\infty}\sum_{j,k=0}^{\infty}\sum \frac{(-i)^{|\alpha|+j+k+\ell}}{\alpha!(j+k+1)!}
{{\rm trace}}_{S(TM)}\Big[\partial_{x_{m}}^{j}\partial_{\xi'}^{\alpha}\partial_{\xi_{m}}^{k}\sigma_{r}^{+}
(fD^{-1}f^{-1}D^{-1})(x',0,\xi',\xi_{m})\nonumber\\
&&\times\partial_{x_{m}}^{\alpha}\partial_{\xi_{m}}^{j+1}\partial_{x_{m}}^{k}\sigma_{l}
\Big((fD^{-1}f^{-1}D^{-1})^{n}\Big)(x',0,\xi',\xi_{m})\Big]
d\xi_{m}\sigma(\xi')dx' ,
\end{eqnarray}
and the sum is taken over $r-k+|\alpha|+\ell-j-1=-(2n+4),r\leq-2,\ell\leq-2n$.

Then, by
Theorem \ref{thm2.3} and direct computations, we have the following theorem.
\begin{thm} \label{thm3.3}
If $M$ is a $2n+4$-dimensional compact oriented manifolds without boundary, then the following equality holds:
\begin{eqnarray}\label{fum3.9}
&&{\rm Wres}\Bigg[\bigg(fD^{-1}f^{-1}D^{-1}\bigg)^{n+1}\Bigg]\nonumber\\
&=&\frac{(2\pi)^{2n+6}}{(2n+4)!}\int_{M}2^{2n+6}{\rm trace}
\bigg\{-\frac{1}{12}s-2f^{-1}\Delta(f)-f^{-2}[|\mathrm{grad}_{M}f|^2+2\Delta(f)]\bigg\}d{\rm Vol_{M}},
\end{eqnarray}
where $s$ is the scalar curvature.
\end{thm}

Locally we can use Theorem \ref{thm3.3} to compute the interior term of (\ref{fum3.7}), then
 \begin{eqnarray}\label{fum3.10}
&&\int_{M}\int_{|\xi|=1}{{\rm trace}}_{S(TM)}\bigg[\sigma_{-n}\Big( (fD^{-1}f^{-1}D^{-1})^{n+1}\Big)\bigg]\sigma(\xi)dx\nonumber\\
&=& \frac{(2\pi)^{2n+6}}{(2n+4)!}\int_{M}2^{2n+6}
\bigg\{-\frac{1}{12}s-2f^{-1}\Delta(f)-f^{-2}\Big[|\mathrm{grad}_{M}f|^2+2\Delta(f)\Big]\bigg\}d{\rm Vol_{M}},
\end{eqnarray}
so we only need to compute $\int_{\partial M}\Phi$.

 When $m=2n+4$ is even, then ${\rm trace}_{S(TM)}[{\rm id}]=2^{\frac{m}{2}}$, the sum is taken over $r-k+|\alpha|+\ell-j=-2n-3,r\leq-2,\ell\leq-2n$, then we have the $\int_{\partial{M}}\Phi$
is the sum of the following five cases:
~\\

\noindent  {\bf case (a)~(I)}~$r=-2,~l=-2n,~k=j=0,~|\alpha|=1$\\

\noindent
By (\ref{fum3.81}), we get
 \begin{eqnarray}\label{case:1}
&&{\rm case~(a)~(I)}\nonumber\\
&=&-\int_{|\xi'|=1}\int^{+\infty}_{-\infty}\sum_{|\alpha|=1}{\rm trace}
\Bigg[\partial^\alpha_{\xi'}\pi^+_{\xi_m}\sigma_{-2}(fD^{-1}f^{-1}D^{-1})\times\partial^\alpha_{x'}
\partial_{\xi_m}\sigma_{-2n}\bigg((fD^{-1}f^{-1}D^{-1})^{n}\bigg)\Bigg](x_0)\nonumber\\
&&\times d\xi_m\sigma(\xi')dx',
\end{eqnarray}
By Lemma \ref{le:32}, for $i<m$, then
\begin{equation}\label{case:2}
\partial_{x_i}\Bigg(\sigma_{-2n}\bigg((fD^{-1}f^{-1}D^{-1})^{n}\bigg)\Bigg)(x_0)
=\partial_{x_i}\bigg(|\xi|^{-2n}\bigg)(x_0)
=(-n)|\xi|^{-2n-2}\partial_{x_i}(|\xi|^2)(x_0)=0,
\end{equation}
\noindent so {\bf case (a)~(I)} vanishes.\\

 \noindent  {\bf case (a)~(II)}~$r=-2,~l=-2n,~k=|\alpha|=0,~j=1$\\

\noindent By (\ref{fum3.81}), we get
\begin{eqnarray}
&&{\rm case~(a)~(II)}\nonumber\\
&=&-\frac{1}{2}\int_{|\xi'|=1}\int^{+\infty}_{-\infty} {\rm trace}
\Bigg[\partial_{x_m}\pi^+_{\xi_m}\sigma_{-2}(fD^{-1}f^{-1}D^{-1})\times
\partial_{\xi_m}^2\sigma_{-2n}\bigg((fD^{-1}f^{-1}D^{-1})^{n}\bigg)\Bigg](x_0)\nonumber\\
&&\times d\xi_m\sigma(\xi')dx'.
\end{eqnarray}
By Lemma \ref{le:32}, we have
\begin{equation}
\partial_{x_m}\sigma_{-2}(fD^{-1}f^{-1}D^{-1})(x_0)|_{|\xi'|=1}=-\frac{h'(0)}{(1+\xi_m^2)^2}.
\end{equation}
By the Cauchy integral formula, then
\begin{eqnarray}
\pi^+_{\xi_m}\partial_{x_m}\sigma_{-2}(fD^{-1}f^{-1}D^{-1})(x_0)|_{|\xi'|=1}
&=&-h'(0)\frac{1}{2\pi i} \lim_{u\rightarrow 0^-}\int_{\Gamma^+}\frac{\frac{1}{(\eta_m+i)^2
(\xi_m+iu-\eta_m)}}{(\eta_m-i)^2}d\eta_m\nonumber\\
&=&\frac{h'(0)(i\xi_m+2)}{4(\xi_m-i)^2}.
\end{eqnarray}
From Lemma \ref{lem3}, we have
\begin{eqnarray}
&&\partial_{\xi_m}^2\sigma_{-2n}\bigg((fD^{-1}f^{-1}D^{-1})^{n}\bigg)(x_0)=\partial^2_{\xi_m}\big(|\xi|^{-2 n}\big) (x_0)\nonumber\\
&=&n(n+1)(|\xi|^{2})^{-n-2} \big(\partial_{\xi_m}|\xi|^2\big)^{2}(x_0)-n(|\xi|^{2})^{-n-1} \partial^2_{\xi_m}\big(|\xi|^2 (x_0)\big)\nonumber\\
&=& \Big((4n+2) \xi_m^{2}-2 \Big)n(1+\xi_m^{2})^{(-n-2)}.
\end{eqnarray}
\noindent We note that
\begin{eqnarray}
&& \int_{-\infty}^{\infty}\bigg[\frac{h'(0)(i\xi_m+2)}{4(\xi_m-i)^2}\times
 \Big((4n+2) \xi_m^{2}-2 \Big)n(1+\xi_m^{2})^{(-n-2)}\bigg]d\xi_m\nonumber\\
&=&\frac{h'(0)\cdot n}{4}\int_{\Gamma^+}\frac{(4n+2)i\xi_m^3+(4+8n)\xi_m^2-2i\xi_m-4}
   {(\xi_m-i)^{(n+4)}(\xi_m+i)^{(n+2)}}d\xi_m\nonumber\\
&=&\frac{h'(0)\cdot n}{4} \frac{2\pi i}{(n+3)!}
\left[\frac{(4n+2)i\xi_m^3+(4+8n)\xi_m^2-2i\xi_m-4}
   {(\xi_m+i)^{(n+2)}}\right]
 ^{(n+3)}|_{\xi_m=i}.
\end{eqnarray}
Since $m=2n+4$ is even, ${\rm trace}_{S(TM)}[{\rm id}]
={\rm dim}(\wedge^*(\frac{2n+4}{2} ))=2^{n+2}.$
Then we obtain
 \begin{equation}
{\rm {\bf case~(a)~(II)}}= \frac{-2^{n}h'(0)\cdot n\pi i}{(n+3)!} {\rm Vol}(S_{2n+2}) \left[\frac{(4n+2)i\xi_m^3+(4+8n)\xi_m^2-2i\xi_m-4}
   {(\xi_m+i)^{(n+2)}}\right]
 ^{(n+3)}|_{\xi_m=i}dx',
\end{equation}
 where ${\rm Vol}(S_{2n+2})$ is the canonical volume of $S_{2n+2}$.\\

\noindent  {\bf case (a)~(III)}~$r=-2,~l=-2n,~j=|\alpha|=0,~k=1$\\

\noindent By (\ref{fum3.81}), we get
\begin{eqnarray}
&&{\rm case~(a)~(III)}\nonumber\\
&=&-\frac{1}{2}\int_{|\xi'|=1}\int^{+\infty}_{-\infty}
{\rm trace}\Bigg [\partial_{\xi_m}\pi^+_{\xi_m}\sigma_{-2}(fD^{-1}f^{-1}D^{-1})\times
\partial_{\xi_m}\partial_{x_m}\sigma_{-2n}\bigg((fD^{-1}f^{-1}D^{-1})^{n}\bigg)\Bigg](x_0)\nonumber\\
&&\times d\xi_m\sigma(\xi')dx'\nonumber\\
&=&\frac{1}{2}\int_{|\xi'|=1}\int^{+\infty}_{-\infty}
{\rm trace}\Bigg [\partial^{2}_{\xi_m}\pi^+_{\xi_m}\sigma_{-2}(fD^{-1}f^{-1}D^{-1})\times
\partial_{x_m}\sigma_{-2n}\bigg((fD^{-1}f^{-1}D^{-1})^{n}\bigg)\Bigg](x_0)\nonumber\\
&&\times d\xi_m\sigma(\xi')dx'.
\end{eqnarray}
\\
By Lemma \ref{lem3}, we have
\begin{equation}\label{fum3.23}
\partial_{\xi_m}^2\pi_{\xi_m}^+\sigma_{-2}(fD^{-1}f^{-1}D^{-1})(x_0)|_{|\xi'|=1}=\frac{-i}{(\xi_m-i)^3},
\end{equation}
and
\begin{equation}\label{fum3.24}
\partial_{x_m} \Bigg(\sigma_{-2n}
\Big((fD^{-1}f^{-1}D^{-1})^{n}\Big)\Bigg)(x_0)
=\partial_{x_m}\bigg((|\xi|^2)^{-n}\bigg) (x_0)
=h'(0)(-n)(1+\xi_{m}^{2})^{-n-1}.
\end{equation}
\\
Then
\begin{eqnarray}\label{fum3.25}
&& \int_{-\infty}^{\infty}{\rm trace}\Big[\frac{-i}{(\xi_m-i)^3}\times h'(0)(-n)(1+\xi_{m}^{2})^{-n-1}\Big] d\xi_m\nonumber\\
&=&i\cdot n\cdot h'(0)\cdot 2^{n+2}
      \int_{\Gamma^+} \frac{1}
   {(\xi_m-i)^{(n+4)}(\xi_m+i)^{(n+1)}}d\xi_m\nonumber\\
&=& -n\cdot h'(0)\cdot 2^{n+3} \cdot\frac{\pi }{(n+3)!}
\left[\frac{1}{(\xi_m+i)^{n+1}}\right]
 ^{(n+3)}|_{\xi_m=i}.
\end{eqnarray}
Then
  \begin{equation}\label{fum3.26}
{\rm {\bf case~(a)~(III)}}=-\frac{2^{n+2}\pi n h'(0)}{(n+3)!} {\rm Vol}(S_{2n+2})\left[\frac{1}{(\xi_m+i)^{n+1}}\right]
 ^{(n+3)}|_{\xi_m=i}dx'.
\end{equation}

\noindent  {\bf case (b)}~$r=-2,~l=-2n-1,~k=j=|\alpha|=0$\\

\noindent By (\ref{fum3.81}) and an integration by parts,, we get
\begin{eqnarray}\label{fum3.27}
&&{\rm case~(b)}\nonumber\\
&=&-i\int_{|\xi'|=1}\int^{+\infty}_{-\infty}{\rm trace} \Bigg[\pi^+_{\xi_m}\sigma_{-2}(fD^{-1}f^{-1}D^{-1})\times
\partial_{\xi_m}\sigma_{-2n-1}\bigg((fD^{-1}f^{-1}D^{-1})^{n}\bigg)\Bigg](x_0)\nonumber\\
&&\times d\xi_m\sigma(\xi')dx'\nonumber\\
&=&i\int_{|\xi'|=1}\int^{+\infty}_{-\infty}{\rm trace} \Bigg[\partial_{\xi_m}\pi^+_{\xi_m}\sigma_{-2}(fD^{-1}f^{-1}D^{-1})\times
\sigma_{-2n-1}\bigg((fD^{-1}f^{-1}D^{-1})^{n}\bigg)\Bigg](x_0)\nonumber\\
&&\times d\xi_m\sigma(\xi')dx'.
\end{eqnarray}
\\
By Lemma \ref{lem3}, we have
\begin{equation}\label{fum3.22}
\partial_{\xi_m}\pi_{\xi_m}^+\sigma_{-2}(fD^{-1}f^{-1}D^{-1})(x_0)|_{|\xi'|=1}=\frac{i}{2(\xi_m-i)^2}.
\end{equation}
Using the recursion formula (4.20) in \cite{KW}, we get
\begin{equation}\label{fum3.28}
\sigma_{3-n}(D^{-n+4})(x,\xi)=\sigma_{5-n}(D^{-n+6})\sigma_{2}^{-1}+\sigma_{2}^{(-\frac{n}{2}+3)}\sigma_{-3}(D^{-2})
             -\sqrt{-1}\partial_{\xi_{\mu}}\sigma_{2}^{(-\frac{n}{2}+3)}\partial_{x_{\mu}}\sigma_{2}^{-1}.
\end{equation}
Then we obtain by induction
\begin{eqnarray}\label{fum3.32}
&&\sigma_{-2n-1}\bigg((fD^{-1}f^{-1}D^{-1})^{n}\bigg)(x,\xi)\nonumber\\
            &=& \sigma_{-2n +1}\bigg((fD^{-1}f^{-1}D^{-1})^{n-1}\bigg)\sigma_{2}^{-1} +\sigma_{2}^{(1-n)}\sigma_{-3}(fD^{-1}f^{-1}D^{-1})
                  -i\cdot\partial_{\xi_{\mu}}\sigma_{2}^{(1-n)}\partial_{x_{\mu}}\sigma_{2}^{-1}\nonumber\\
          &=& n\cdot\sigma_{2}^{(1-n)}\sigma_{-3}(fD^{-1}f^{-1}D^{-1})-
          i\cdot \sum_{k=0}^{n-2}\partial_{\xi_{\mu}}
\sigma_{2}^{(1-n+k)}\partial_{x_{\mu}}\sigma_{2}^{-1}\big(\sigma_{2}^{-1}\big)^{k}.
\end{eqnarray}

In the normal coordinate, $g^{ij}(x_0)=\delta_i^j$ and $\partial_{x_j}(g^{\alpha\beta})(x_0)=0,$ {\rm if
}$j<m;\partial_{x_j}(g^{\alpha\beta})(x_0)=h'(0)\delta^\alpha_\beta,~{\rm if }~j=m.$ So by Lemma A.2 in \cite{Wa3}, we have $\Gamma^m(x_0)=\frac{2n+3}{2}h'(0)$ and
$\Gamma^k(x_0)=0$ for $k<m$. By the definition of $\delta^k$ and Lemma \ref{lem3}, we have $\delta^m(x_0)=0$ and
$\delta^k=\frac{1}{4}h'(0)
c(\widetilde{e_k})c(\widetilde{e_m})$ for
$k<m$. So
 \begin{eqnarray}\label{fum3.28}
&&\sigma_{-3}(fD^{-1}f^{-1}D^{-1})(x_0)|_{|\xi'|=1}\nonumber\\
&=&-i|\xi|^{-4}\xi_k(\Gamma^k-2\delta^k)(x_0)|_{|\xi'|=1}-i|\xi|^{-6}2\xi^j\xi_\alpha\xi_\beta
\partial_jg^{\alpha\beta}(x_0)|_{|\xi'|=1}-i|\xi|^{-4}\xi_k\big[c(\partial^j)\cdot f\cdot c(df^{-1})\big]\nonumber\\
&=&\frac{i}{(1+\xi_m^2)^2}\Big(\frac{1}{2}h'(0)
\sum_{k<m}\xi_k c(\widetilde{e_k})c(\widetilde{e_n})
-\frac{2n+3}{2}h'(0)\xi_m\Big)-\frac{2ih'(0)\xi_m}{(1+\xi_m^2)^3}
-\frac{i}{(1+\xi_m^2)^2}\sum^{m}_{k=1}\xi_k\big[c(\widetilde{e_k})\nonumber\\
&&\times f\cdot c(df^{-1})\big].
\end{eqnarray}
We note that $\int_{|\xi'|=1}\xi_1\cdots\xi_{2q+1}\sigma(\xi')=0$, so the first term and the fourth term in (\ref{fum3.28}) has no contribution for computing {\bf case (b)}.

On the other hand, we have
\begin{equation}\label{fum3.29}
\sigma_{2}^{(1-n)} (x_0)= (1+\xi_{m}^{2})^{(1-n)},
\end{equation}
and
\begin{equation}\label{fum3.30}
\partial_{x_{j}}(|\xi|^{-2})(x_0)=0,~~~j<m.
\end{equation}
Then
\begin{eqnarray}\label{fum3.31}
&& -i \sum_{k=0}^{n-2}\partial_{\xi_{\mu}}\sigma_{2}^{(1-n+k)}
\partial_{x_{\mu}}\sigma_{2}^{-1}\big(\sigma_{2}^{-1}\big)^{k}(x_0)\nonumber\\
&=&-i\sum_{k=0}^{n-2} \partial_{\xi_{n}}\Big[(|\xi|^{2})^{(1-n+k)}\Big]
\partial_{x_{n}}(|\xi|^{2})^{-1}(|\xi|)^{-2k}=i\sum_{k=0}^{n-2} (|\xi|^{2})^{(-n+k)} (1-n+k)2\xi_{n}|\xi|^{-4}
h'(0) (1+\xi_{m}^{2})^{-k}\nonumber\\
&=&i\sum_{k=0}^{n-2} (1+\xi_{m}^{2})^{(-n)}(1+\xi_{m}^{2})^{-2}
   h'(0)\xi_{n}(-2n+2k+2)=i\sum_{k=0}^{n-2}  h'(0)(-2n+2k+2)\xi_{n}(1+\xi_{m}^{2})^{(-n)}
\nonumber\\
&=&i h'(0)(-n^2+n)  \xi_{n} (1+\xi_{m}^{2})^{(-n-2)}.
\end{eqnarray}

In conclusion, we obtain
\begin{eqnarray}\label{fum3.32}
&&\sigma_{-2n-1}\bigg((fD^{-1}f^{-1}D^{-1})^{n}\bigg)(x,\xi)\nonumber\\
            &=& n(1+\xi_{m}^{2})^{1-n}
   \Big(\frac{-i (2n+3)h'(0)\xi_m}{2(1+\xi_m^2)^2}-\frac{2ih'(0)\xi_m}{(1+\xi_m^2)^3} \Big)-i h'(0)(n^2-n) \xi_{n} (1+\xi_{m}^{2})^{(-n-2)}.
\end{eqnarray}

From (\ref{fum3.22}) and (\ref{fum3.32}), we obtain
\begin{eqnarray}
{\bf case~ (b)}&=& i\int_{|\xi'|=1}\int^{+\infty}_{-\infty}{\rm trace}
\Bigg\{\frac{i}{2(\xi_m-i)^2}\times \Bigg[n(1+\xi_{m}^{2})^{1-n}
   \Big(\frac{-i}{(1+\xi_m^2)^2}\times\frac{2n+3}{2}h'(0)\xi_m\nonumber\\
   &&-\frac{2ih'(0)\xi_m}{(1+\xi_m^2)^3} \Big)-i h'(0)(n^2-n) \xi_{n} (1+\xi_{m}^{2})^{(-n-2)}\Bigg]
\Bigg\}d\xi_m\sigma(\xi')dx' \nonumber\\
&=&-\frac{i h'(0)}{8} {\rm Vol}(S_{2n+2})\int_{\Gamma^+}\frac{-2n(2n+3)\xi_{n}^{3}+(-8n^2-10n-4)  \xi_{n}}{(\xi_m-i)^{n+4}(\xi_m+i)^{n+2}}d\xi_ndx'\nonumber\\
&=&-\frac{i h'(0)}{8} {\rm Vol}(S_{2n+2})2^{n+2}  \frac{2\pi i}{(n+3)!}
  \left[\frac{-2n(2n+3)\xi_{n}^{3}-(8n^2+10n+4)\xi_{n}}
  {(\xi_m+i)^{n+2}}\right]^{n+3}|_{\xi_m=i}dx'\nonumber\\
&=&\frac{2^{n}\pi h'(0)}{(n+3)!} {\rm Vol}(S_{2n+2})
  \left[\frac{-2n(2n+3)\xi_{n}^{3}-(8n^2+10n+4)\xi_{n}}
  {(\xi_m+i)^{n+2}}\right]^{(n+3)}|_{\xi_m=i}dx'.
\end{eqnarray}

\noindent {\bf  case (c)}~$r=-3,~l=-2n,~k=j=|\alpha|=0$\\

By (\ref{fum3.81}), we get
\begin{eqnarray}
{\rm case~ (c)}&=&-i\int_{|\xi'|=1}\int^{+\infty}_{-\infty}{\rm trace} \Bigg[\pi^+_{\xi_m}\sigma_{-3}(fD^{-1}f^{-1}D^{-1})\times
\partial_{\xi_m}\sigma_{-2n}\bigg((fD^{-1}f^{-1}D^{-1})^{n}\bigg)\Bigg](x_0)\nonumber\\
&&\times d\xi_m\sigma(\xi')dx'
\end{eqnarray}
By Lemma \ref{lem3}, we have
\begin{equation}\label{fum3.37}
\partial_{\xi_m}\Bigg[\sigma_{-2n}\bigg((fD^{-1}f^{-1}D^{-1})^{n}\bigg)\Bigg](x_0)
=\partial_{\xi_m}\bigg((|\xi|^2)^{-n}\bigg)(x_0)=-2n\xi_m(1+\xi_m^2)^{-n-1}.
\end{equation}
By the Cauchy integral formula, we obtain
\begin{eqnarray}\label{fum3.38}
\pi^+_{\xi_m} \Big( \frac{\xi_m}{(1+\xi_m^2)^{2}}\Big)
&=&\frac{1}{2\pi i} \int_{\Gamma^+}\frac{\eta_m}{(\xi_m-\eta_m)
(1+\eta_m^{2} )^2}d\eta_m\nonumber\\
&=&\left[\frac{\eta_m}{(\xi_m-\eta_m)(\eta_m^{2}+i)^{2}}\right]^{1}\Big|_{\eta_n=i}
     \nonumber\\
&=& \frac{-i}{4(\xi_m-i)^2},
\end{eqnarray}
and
\begin{eqnarray}\label{fum3.39}
\pi^+_{\xi_m} \Big( \frac{\xi_m}{(1+\xi_m^2)^{3}}\Big)
= \frac{-i}{16(\xi_m-i)^2}- \frac{1}{8(\xi_m-i)^3}.
\end{eqnarray}
In conclusion, we obtain
\begin{eqnarray}\label{fum3.40}
\pi^+_{\xi_m}\Big(  \sigma_{-3}(fD^{-1}f^{-1}D^{-1})(x_0)|_{|\xi'|=1}\Big)
             &=& -i h'(0)\pi^+_{\xi_m}\Big(
             \frac{(2n+3)\xi_m }{2(1+\xi_m^2)^2}+ \frac{2\xi_m}{(1+\xi_m^2)^3}\Big)\nonumber\\
         &=& -i h'(0)\Big[(2n+3)\pi^+_{\xi_m}
         \Big(\frac{\xi_m }{(1+\xi_m^2)^2}\Big)+2 \pi^+_{\xi_m} \Big(\frac{\xi_m}{(1+\xi_m^2)^3}\Big) \Big]\nonumber\\
         &=&i h'(0)\Big[\frac{i(2n+2)}{8(\xi_m-i)^2}+ \frac{1}{4(\xi_m-i)^3} \Big].
\end{eqnarray}
Therefore, by (\ref{fum3.37}) and (\ref{fum3.40}), we have
\begin{eqnarray}\label{fum3.41}
{\bf case~ (c)}&=& -i\int_{|\xi'|=1}\int^{+\infty}_{-\infty}{\rm trace}
\Bigg[ \frac{-2n\xi_m}{(1+\xi_m^2)^{n+1}}
\times i h'(0)\Big[\frac{i(2n+2)}{8(\xi_m-i)^2}+ \frac{1}{4(\xi_m-i)^3}\Big]
\Bigg]d\xi_m\sigma(\xi')dx' \nonumber\\
&=& (-n)2^{n} {\rm Vol}(S_{2n+2})h'(0)
  \int_{\Gamma^+}\frac{\big(2i(n+1)\xi_{m}+2n+4\big)\xi_{n}}
  {(\xi_m+i)^{n+1}(\xi_m-i)^{n+4}}d\xi_mdx' \nonumber\\
&=&(-n)2^{n} {\rm Vol}(S_{2n+2})h'(0) \frac{2\pi i}{(n+3)!}
  \left[\frac{\big(2i(n+1)\xi_{n}+2n+4\big)\xi_{m}}
  {(\xi_m+i)^{n+1}}\right]
  ^{(n+3)}|_{\xi_m=i}dx'.
\end{eqnarray}
Since $\Phi$ is the sum of the {\bf case~(a)}, {\bf case~(b)} and {\bf case~(c)}, so
\begin{eqnarray}\label{fum3.42}
\Phi&=&\frac{2^{n}h'(0)\pi ni}{(n+3)!}{\rm Vol}(S_{2n+2})dx'\cdot
\left[\frac{-4ni \xi_m^{3}-4(2n+2)\xi_m^{2}+(4ni+8i)\xi_m}
  {(\xi_m+i)^{n+2}}\right]
  ^{(n+3)}|_{\xi_m=i}\nonumber\\
 &:=&\frac{2^{n}h'(0)\pi ni}{(n+3)!}{\rm Vol}(S_{2n+2})dx'\cdot Q_{0},
\end{eqnarray}
where
\begin{eqnarray}\label{fum9.1}
Q_0&=&\left[\frac{-4ni \xi_m^{3}-4(2n+2)\xi_m^{2}+(4ni+8i)\xi_m}
  {(\xi_m+i)^{n+2}}\right]
  ^{(n+3)}|_{\xi_m=i}\nonumber\\
&=&(-1)^{n}(n+3)!(1+i)4^{-(n+1)}
\bigg[(2+2i)nC_{-n-2}^{n}+(1+i)(n-2)C_{-n-2}^{n+1}
-(3+i)C_{-n-2}^{n+2}\nonumber\\
&&-C_{-n-2}^{n+3}\bigg].
\end{eqnarray}
Combining (\ref{fum3.10}) and (\ref{fum3.42}), we obtain Theorem 1.1.

\section{The noncommutative residue $\widetilde{{\rm Wres}}\bigg[\pi^+(fD^{-1})\circ\pi^+\Big((f^{-1}D^{-1})\cdot(fD^{-1}f^{-1}D^{-1})^{n}\Big)\bigg]$ on even dimensional manifolds with boundary}
\label{section:4}
In the following, we will compute the residue $\widetilde{{\rm Wres}}\bigg[\pi^+(fD^{-1})\circ\pi^+\Big((f^{-1}D^{-1})\cdot(fD^{-1}f^{-1}D^{-1})^{n}\Big)\bigg]$ for nonzero
smooth functions $f,~f^{-1}$ on even dimensional oriented
compact spin manifolds with boundary and get a general Kastler-Kalau-Walze
type theorem in this case. By Theorem \ref{thm2.1}, we have
\begin{eqnarray}\label{fum4.1.}
&&\widetilde{{\rm Wres}}\bigg[\pi^+(fD^{-1})\circ\pi^+\Big((f^{-1}D^{-1})\cdot(fD^{-1}f^{-1}D^{-1})^{n}\Big)\bigg]\nonumber\\
&=&\int_{M}\int_{|\xi|=1}{{\rm trace}}_{S(TM)}\big[\sigma_{-n}\big( (fD^{-1}f^{-1}D^{-1})^{n+1}\big)\big]\sigma(\xi)dx+\int_{\partial M}\Psi,
\end{eqnarray}
where
\begin{eqnarray}\label{fum4.2}
\Psi&=&\int_{|\xi'|=1}\int_{-\infty}^{+\infty}\sum_{j,k=0}^{\infty}\sum \frac{(-i)^{|\alpha|+j+k+\ell}}{\alpha!(j+k+1)!}
{{\rm trace}}_{S(TM)}\Big[\partial_{x_{m}}^{j}\partial_{\xi'}^{\alpha}\partial_{\xi_{m}}^{k}\sigma_{r}^{+}
(fD^{-1})(x',0,\xi',\xi_{m})\nonumber\\
&&\times\partial_{x_{m}}^{\alpha}\partial_{\xi_{m}}^{j+1}\partial_{x_{m}}^{k}\sigma_{l}
\Big((f^{-1}D^{-1})\cdot
(fD^{-1}f^{-1}D^{-1})^{n}\Big)(x',0,\xi',\xi_{m})\Big]
d\xi_{m}\sigma(\xi')dx' ,
\end{eqnarray}
and the sum is taken over $r-k+|\alpha|+\ell-j-1=-(2n+4),r\leq-1,\ell\leq-2n-1$.

Locally we can use Theorem \ref{thm3.3} to compute the interior term of (\ref{fum3.7}), then
 \begin{eqnarray}\label{fum4.3}
&&\int_{M}\int_{|\xi|=1}{{\rm trace}}_{S(TM)}\bigg[\sigma_{-n}\Big( (fD^{-1}f^{-1}D^{-1})^{n+1}\Big)\bigg]\sigma(\xi)dx\nonumber\\
&=& \frac{(2\pi)^{2n+6}}{(2n+4)!}\int_{M}2^{2n+6}
\bigg\{-\frac{1}{12}s-2f^{-1}\Delta(f)-f^{-2}\Big[|\mathrm{grad}_{M}f|^2+2\Delta(f)\Big]\bigg\}d{\rm Vol_{M}},
\end{eqnarray}
so we only need to compute $\int_{\partial M}\Psi$.

When $m=2n+4$ is even, then ${\rm trace}_{S(TM)}[{\rm id}]=2^{\frac{m}{2}}$, the sum is taken over $r-k+|\alpha|+\ell-j=-2n-3,r\leq-1,\ell\leq-2n-1$, then we have the $\int_{\partial{M}}\Psi$
is the sum of the following five cases:
~\\

\noindent  {\bf case ~(1)}~$r=-1, l=-2n-1, j=k=0, |\alpha|=1$.
~\\

\noindent By (\ref{fum4.2}), we get
 \begin{eqnarray}\label{fum4.4}
&&{\rm case~(1)}\nonumber\\
&=&-\int_{|\xi'|=1}\int^{+\infty}_{-\infty}\sum_{|\alpha|=1}{\rm trace}
\Bigg[\partial^{\alpha}_{\xi'}\pi^{+}_{\xi_{m}}\sigma_{-1}(fD^{-1})
      \times\partial^{\alpha}_{x'}\partial_{\xi_{m}}\sigma_{-2n-1}
      \bigg((f^{-1}D^{-1})\cdot (fD^{-1}\cdot f^{-1}D^{-1})^{n}\bigg)\Bigg](x_0)\nonumber\\
      &&\times d\xi_m\sigma(\xi')dx'\nonumber\\
&=&-\int_{|\xi'|=1}\int^{+\infty}_{-\infty}\sum_{|\alpha|=1}{\rm trace}
\Bigg[\partial^{\alpha}_{\xi'}\pi^{+}_{\xi_{m}}\big(f\sigma_{-1}(D^{-1})\big)
      \times\partial^{\alpha}_{x'}\partial_{\xi_{m}}\bigg(f^{-1}\sigma_{-2n-1}
(D^{-2n-1})\bigg)\Bigg](x_0)d\xi_m\sigma(\xi')dx'\nonumber\\
&=&-\int_{|\xi'|=1}\int^{+\infty}_{-\infty}\sum_{|\alpha|=1}{\rm trace}
\Bigg\{f\cdot\partial^{\alpha}_{\xi'}\pi^{+}_{\xi_{m}}\big(\sigma_{-1}(D^{-1})\big)
      \times\bigg[\big(\partial^{\alpha}_{x'}(f^{-1})\big)\partial_{\xi_{m}}\bigg(\sigma_{-2n-1}
(D^{-2n-1})\bigg)+f^{-1}\nonumber\\
&&\times\partial^{\alpha}_{x'}\partial_{\xi_{m}}\bigg(\sigma_{-2n-1}
(D^{-2n-1})\bigg)\bigg]\Bigg\}(x_0)d\xi_m\sigma(\xi')dx'\nonumber\\
&=&-\int_{|\xi'|=1}\int^{+\infty}_{-\infty}\sum_{|\alpha|=1}{\rm trace}
\Bigg[\partial^{\alpha}_{\xi'}\pi^{+}_{\xi_{m}}\big(\sigma_{-1}(D^{-1})\big)
      \times\partial^{\alpha}_{x'}\partial_{\xi_{m}}\bigg(\sigma_{-2n-1}(D^{-2n-1})\bigg)\Bigg](x_0)d\xi_m\sigma(\xi')dx'
\nonumber\\
&&-f\cdot\sum_{j<1}\partial_{x_{j}}(f^{-1})\int_{|\xi'|=1}\int^{+\infty}_{-\infty}\sum_{|\alpha|=1}{\rm trace}\Bigg[\partial^{\alpha}_{\xi'}\pi^{+}_{\xi_{m}}\big(\sigma_{-1}(D^{-1})\big)
      \times\partial_{\xi_{m}}\bigg(\sigma_{-2n-1}(D^{-2n-1})\bigg)\Bigg](x_0)\nonumber\\
&&\times d\xi_m\sigma(\xi')dx'.
\end{eqnarray}
By Lemma 2.2 in \cite{Wa3} and (3.12) in \cite{WJ5}, we have for $j<m$
\begin{eqnarray}\label{fum4.5}
&&\partial_{x_j}\sigma_{-2n-1}(D^{-2n-1})(x_0)=\partial_{x_j}\big[\sqrt{-1}c(\xi)|\xi|^{-2n-2}\big]\nonumber\\
&=&\sqrt{-1}\big[\partial_{x_j}c(\xi)\big](x_0)|\xi|^{-2n-2}
+\sqrt{-1}c(\xi)\partial_{x_j}(|\xi|^{-2n-2})(x_0)=0,
\end{eqnarray}
 so
\begin{eqnarray}\label{fum4.6}
-\int_{|\xi'|=1}\int^{+\infty}_{-\infty}\sum_{|\alpha|=1}{\rm trace}
\Bigg[\partial^{\alpha}_{\xi'}\pi^{+}_{\xi_{m}}\big(\sigma_{-1}(D^{-1})\big)
      \times\partial^{\alpha}_{x'}\partial_{\xi_{m}}\bigg(\sigma_{-2n-1}(D^{-2n-1})\bigg)\Bigg]
      (x_0)d\xi_m\sigma(\xi')dx'
=0.
\end{eqnarray}

\noindent By Lemma \ref{lem3} and direct calculations, for $i<m$, we obtain
\begin{eqnarray}\label{fum4.7}
&&\partial^{\alpha}_{\xi'}\pi^{+}_{\xi_{m}}\sigma_{-1}(D^{-1})(x_0)|_{|\xi'|=1}
=\partial_{\xi_i}\pi^{+}_{\xi_{m}}\sigma_{-1}(D^{-1})(x_0)|_{|\xi'|=1}\nonumber\\
&=&\frac{c(dx_i)}{2(\xi_{m}-i)}-\frac{\xi_i(\xi_{m}-2i)c(\xi')+\xi_ic(dx_m)}{2(\xi_{m}-i)^2},
\end{eqnarray}
and we get
\begin{eqnarray}\label{fum4.8}
\partial_{\xi_{m}}\bigg(\sigma_{-2n-1}(D^{-2n-1})\bigg)
=\frac{\sqrt{-1}c(dx_m)}{|\xi|^{2n+2}}-\frac{(2n+2)\sqrt{-1}\Big[\xi_nc(\xi')+\xi^2_nc(dx_m)\Big]}{|\xi|^{2n+4}}.
\end{eqnarray}
Then for $i<m$, we have
\begin{eqnarray}\label{fum4.9}
&&{\rm trace}
\Bigg[\partial^{\alpha}_{\xi'}\pi^{+}_{\xi_{m}}\sigma_{-1}(D^{-1})
      \times\partial_{\xi_{m}}\bigg(\sigma_{-2n-1}(D^{-2n-1})\bigg)\Bigg](x_0)\nonumber\\
&=&-\xi_i{\rm trace}
\Big[\frac{c(dx_m)^{2}}{2(\xi_{m}-i)^2|\xi|^{2n+2}}\Big]-4i\xi_{m}\xi_i{\rm trace}
\Big[\frac{c(dx_i)^{2}}{2(\xi_{m}-i)|\xi|^{2n+4}}\Big]+4i\xi_{m}\xi_i(\xi_{m}-2i)\nonumber\\
&&\times{\rm trace}
\Big[\frac{c(\xi')^{2}}{2(\xi_{m}-i)^2|\xi|^{2n+4}}\Big]+4i\xi^{2}_m\xi_i{\rm trace}
\Big[\frac{c(dx_m)^{2}}{2(\xi_{m}-i)^2|\xi|^{2n+4}}\Big].
\end{eqnarray}
We note that $i<m,~\int_{|\xi'|=1}\xi_i\sigma(\xi')=0$,
so
\begin{eqnarray}\label{fum4.10}
&&-f\sum\limits_{j<m}\partial_{x_j}(f^{-1})\int_{|\xi'|=1}\int^{+\infty}_{-\infty}\sum_{|\alpha|=1}{\rm trace}
\Bigg[\partial^{\alpha}_{\xi'}\pi^{+}_{\xi_{m}}\sigma_{-1}(D^{-1})
      \times\partial_{\xi_{m}}\bigg(\sigma_{-2n-1}(D^{-2n-1})\bigg)\Bigg](x_0)\nonumber\\
      &&\times d\xi_{m}\sigma(\xi')dx'\nonumber\\
&=&0.
\end{eqnarray}
Then we have ${\bf case~(1)}=0$.
~\\

\noindent  {\bf case ~(2)}~$r=-1, l=-2n-1, |\alpha|=k=0, j=1$.\\

\noindent By (\ref{fum4.2}), we have
\begin{eqnarray}\label{fum4.11}
&&{\rm case~(2)}\nonumber\\
&=&-\frac{1}{2}\int_{|\xi'|=1}\int^{+\infty}_{-\infty} {\rm
trace} \Bigg[\partial_{x_{m}}\pi^{+}_{\xi_{m}}\sigma_{-1}(fD^{-1})
      \times\partial^{2}_{\xi_{m}}\sigma_{-2n-1}
      \bigg((f^{-1}D^{-1})\cdot (fD^{-1}\cdot f^{-1}D^{-1})^{n}\bigg)\Bigg](x_0)\nonumber\\
      &&\times d\xi_{m}\sigma(\xi')dx'\nonumber\\
&=&-\frac{1}{2}\int_{|\xi'|=1}\int^{+\infty}_{-\infty}{\rm trace}
\Bigg[\partial_{x_{m}}\pi^{+}_{\xi_{m}}\sigma_{-1}(D^{-1})
      \times\partial^2_{\xi_{m}}\sigma_{-2n-1}
      \bigg(D^{-2n-1}\bigg)\Bigg](x_0)d\xi_{m}\sigma(\xi')dx'\nonumber\\
      &&-\frac{1}{2}f^{-1}\partial_{x_{m}}(f)
      \int_{|\xi'|=1}\int^{+\infty}_{-\infty}{\rm trace}
\Bigg[\pi^{+}_{\xi_{m}}\sigma_{-1}(D^{-1})
      \times\partial^2_{\xi_{m}}\sigma_{-2n-1}
\bigg(D^{-2n-1}\bigg)\Bigg](x_0)d\xi_{m}\sigma(\xi')dx'\nonumber\\
&=&-\frac{1}{2}\int_{|\xi'|=1}\int^{+\infty}_{-\infty}{\rm trace}
\Bigg[\partial^2_{\xi_{m}}\partial_{x_{m}}\pi^{+}_{\xi_{m}}\sigma_{-1}(D^{-1})
      \times\sigma_{-2n-1}
      \bigg(D^{-2n-1}\bigg)\Bigg](x_0)d\xi_{m}\sigma(\xi')dx'-\frac{1}{2}\nonumber\\
      &&\times f^{-1}\partial_{x_{m}}(f)
      \int_{|\xi'|=1}\int^{+\infty}_{-\infty}{\rm trace}
\Bigg[\pi^{+}_{\xi_{m}}\sigma_{-1}(D^{-1})
      \times\partial^2_{\xi_{m}}\sigma_{-2n-1}
\bigg(D^{-2n-1}\bigg)\Bigg](x_0)d\xi_{m}\sigma(\xi')dx'.
\end{eqnarray}
By (2.2.23) in \cite{Wa3}, we have
\begin{eqnarray}\label{fum4.12}
&&\pi^+_{\xi_{m}}\partial_{x_m}\sigma_{-1}(D^{-1})(x_0)|_{|\xi'|=1}\nonumber\\
&=&\frac{\partial_{x_m}\big[c(\xi')\big](x_0)}{2(\xi_{m}-i)}+\sqrt{-1}h'(0)
\left[\frac{ic(\xi')}{4(\xi_{m}-i)}+\frac{c(\xi')+ic(dx_m)}{4(\xi_{m}-i)^2}\right].
\end{eqnarray}
So
\begin{eqnarray}\label{fum4.13}
&&\partial_{\xi_{m}}^2\pi^+_{\xi_{m}}\partial_{x_m}\sigma_{-1}(D^{-1})(x_0)|_{|\xi'|=1}\nonumber\\
&=&\frac{\partial_{x_m}\big[c(\xi')\big](x_0)}{4(\xi_{m}-i)^3}+\sqrt{-1}h'(0)
\left[\frac{ic(\xi')}{8(\xi_{m}-i)^3}+\frac{c(\xi')+ic(dx_m)}{24(\xi_{m}-i)^4}\right].
\end{eqnarray}
\noindent We know that
\begin{eqnarray}\label{fum4.14}
\sigma_{-2n-1}(D^{-2n-1})=\frac{\sqrt{-1}[c(\xi')+\xi_nc(dx_m)]}{(1+\xi_{m}^2)^{n+1}},
\end{eqnarray}
\noindent By the relation of the Clifford action and ${\rm trace}{AB}={\rm trace }{BA}$, then we have the equalities:\\
$${\rm trace}\big[c(\xi')c(dx_m)\big]=0;~~{\rm trace}\big[c(dx_m)^2\big]=-2^{n+2};~~{\rm trace}\big[c(\xi')^2\big](x_0)|_{|\xi'|=1}=-2^{n+2};~~$$
$${\rm trace}\big[\partial_{x_m}c(\xi')c(dx_m)\big]=0;~~{\rm trace}\big[\partial_{x_m}c(\xi')c(\xi')\big](x_0)|_{|\xi'|=1}=-2^{n+1}h'(0).$$
By (\ref{fum4.11}), (\ref{fum4.13}) and (\ref{fum4.14}), we have
\begin{eqnarray}\label{fum4.15}
&& {\rm trace} \Big[ \partial_{\xi_{m}}^2\partial_{x_m}\pi^+_{\xi_{m}}\sigma_{-1}(D^{-1})
\times \sigma_{-2n-1}(D^{-2n-1})\Big](x_0)|_{|\xi'|=1}\nonumber\\
&=& {\rm trace}\bigg\{\left[  \frac{\partial_{x_m}[c(\xi')](x_0)}{4(\xi_{m}-i)^3}+ih'(0)
\Big[\frac{ic(\xi')}{8(\xi_{m}-i)^3}+\frac{c(\xi')+ic(dx_m)}{24(\xi_{m}-i)^4}\Big]\right]\times \frac{i[c(\xi')+\xi_mc(dx_m)]}{(1+\xi_{m}^2)^{n+1}}\bigg\}\nonumber\\
&&\nonumber\\
&=&{\rm trace}\bigg\{\Big[\frac{\partial_{x_m}[c(\xi')](x_0)}{4(\xi_{m}-i)^3}+
\frac{(4i-3\xi_{m})h'(0)}{24(\xi_{m}-i)^4}c(\xi')-\frac{h'(0)c(dx_m)}{24(\xi_{m}-i)^4}\Big]
\times \frac{\sqrt{-1}[c(\xi')+\xi_mc(dx_m)]}{(1+\xi_{m}^2)^{n+1}}\bigg\}\nonumber\\
&=&\frac{2^{n+1}h'(0)i}{12(\xi_{m}+i)^{n+1}(\xi_{m}-i)^{n+4}}.
\end{eqnarray}
Thus, we have
\begin{eqnarray}\label{fum4.16}
&&-\frac{1}{2}\int_{|\xi'|=1}\int^{+\infty}_{-\infty}{\rm trace}
\Bigg[\partial^2_{\xi_{m}}\partial_{x_{m}}\pi^{+}_{\xi_{m}}\sigma_{-1}(D^{-1})
      \times\sigma_{-2n-1}
      \bigg(D^{-2n-1}\bigg)\Bigg](x_0)d\xi_{m}\sigma(\xi')dx'\nonumber\\
&=&{\rm Vol}(S_{2n+2}) \frac{2^{n+1}h'(0)\pi}{12(n+3)!}\Big[(\xi_{m}+i)^{-n-1}\Big]^{(n+3)}|_{\xi_{m}=i}dx',
\end{eqnarray}
where ${\rm Vol}(S_{2n+2})$ is the canonical volume of $S_{2n+2}$ and denote the $p$-th derivative of $f(\xi_{m})$ by $[f(\xi_{m})]^{(p)}$.\\
By (4.14) and direct calculations, we have\\
\begin{equation}\label{fum4.17}
\partial_{\xi_{m}}\sigma_{-2n-1}(D^{-2n-1})
=\frac{-2(n+1)i\xi_{m}[c(\xi')+\xi_nc(dx_m)]}{(1+\xi_{m}^{2})^{n+2}}
\end{equation}
and
\begin{eqnarray}\label{fum4.18}
&&\partial^{2}_{\xi_{m}}\sigma_{-2n-1}(D^{-2n-1})\nonumber\\
&=&i\bigg[\frac{4\xi^{2}_{m}(n+1)(n+2)\big(c(\xi')+\xi_{m}c(dx_{m})\big)}{(1+\xi_{m}^{2})^{n+3}}
-\frac{ 6\xi_{m}c(dx_{m})(n+1)+ 2(n+1)c(\xi') }{ (1+\xi_{m}^{2})^{n+2} }\bigg].
\end{eqnarray}

On the other hand, by calculations, we have
\begin{equation}\label{fum4.19}
\pi^{+}_{\xi_{m}}\sigma_{-1}(D^{-1})(x_{0})|_{|\xi'|=1}
=-\frac{c(\xi')+\sqrt{-1}c(dx_{m})}{2(\xi_{m}-\sqrt{-1})}.
\end{equation}
By (\ref{fum4.11}), (\ref{fum4.18}) and (\ref{fum4.19}), we get
\begin{eqnarray}\label{fum4.20}
&&{\rm
trace} \Big[\pi^{+}_{\xi_{m}}\sigma_{-1}(D^{-1})
      \times\partial^{2}_{\xi_{m}}\sigma_{-2n-1}(D^{-2n-1})\Big](x_0)\nonumber\\
&=&\frac{2^{n+2}(n+1)\Big[2i\xi^{2}_m(n+2)(1+\xi_{m})-(2\xi_{m}-i-i\xi_{m})(1+\xi^{2}_m)\Big]}
{(\xi_{m}-i)^{n+4}(\xi_{m}+i)^{n+3}},
\end{eqnarray}
then we have
\begin{eqnarray}\label{fum4.21}
&&-\frac{1}{2}f^{-1}\partial_{x_{m}}(f)
      \int_{|\xi'|=1}\int^{+\infty}_{-\infty}{\rm trace}
\Bigg[\pi^{+}_{\xi_{m}}\sigma_{-1}(D^{-1})
      \times\partial^2_{\xi_{m}}\sigma_{-2n-1}
\bigg(D^{-2n-1}\bigg)\Bigg](x_0)d\xi_{m}\sigma(\xi')dx'\nonumber\\
&=&\frac{-\pi i f^{-1}\partial_{x_{m}}(f)}{(n+3)!}{\rm Vol}(S_{2n+2})
G_0dx',
\end{eqnarray}
where $$G_0=\bigg[\frac{2^{n+2}(n+1)\big[(2ni+3i+2)\xi^{3}_m+(2n+3)i\xi^{2}_m+(2-i)\xi_{m}-i\big]}{(\xi_{m}+1)^{n+3}}  \bigg]^{(n+3)}|_{\xi_{m}=i}.$$
\\
Combining (\ref{fum4.11}), (\ref{fum4.16}) and (\ref{fum4.21}), we obtain
\begin{eqnarray*}
{\bf case~(2)}&=&\bigg[ \frac{2^{n+1}h'(0)\pi}{12(n+3)!}\Big(\frac{1}{(\xi_{m}+i)^{n+1}}\Big)^{(n+3)}|_{\xi_{m}=i}-
\frac{\pi i f^{-1}\partial_{x_{m}}(f)}{(n+3)!}
G_0\bigg]{\rm Vol}(S_{2n+2})dx'.
\end{eqnarray*}

\noindent  {\bf case~(3)}~$r=-1,l=-2n-1,|\alpha|=j=0,k=1$.\\

\noindent By (\ref{fum4.2}), we have
 \begin{eqnarray}\label{fum4.22}
&&{\rm case~ (a)~(3)}\nonumber\\
&=&-\frac{1}{2}\int_{|\xi'|=1}\int^{+\infty}_{-\infty}{\rm trace} \Big[\partial_{\xi_{m}}\pi^{+}_{\xi_{m}}\sigma_{-1}(fD^{-1})
      \times\partial_{\xi_{m}}\partial_{x_{m}}
      \sigma_{-2n-1}\bigg((f^{-1}D^{-1})\cdot (fD^{-1}\cdot f^{-1}D^{-1})^{n}\bigg)\Big](x_0)\nonumber\\
      &&\times d\xi_{m}\sigma(\xi')dx'\nonumber\\
&=&-\frac{1}{2}\int_{|\xi'|=1}\int^{+\infty}_{-\infty}{\rm trace}
\Bigg[\partial_{\xi_{m}}\pi^{+}_{\xi_{m}}\sigma_{-1}(D^{-1})
      \times\partial_{\xi_{m}}\partial_{x_{m}}\sigma_{-2n-1}
      \bigg(D^{-2n-1}\bigg)\Bigg](x_0)d\xi_{m}\sigma(\xi')dx'\nonumber\\
      &&-\frac{1}{2}f\cdot\partial_{x_{m}}(f^{-1})
      \int_{|\xi'|=1}\int^{+\infty}_{-\infty}{\rm trace}
\Bigg[\partial_{\xi_{m}}\pi^{+}_{\xi_{m}}\sigma_{-1}(D^{-1})
      \times\partial_{\xi_{m}}\sigma_{-2n-1}
\bigg(D^{-2n-1}\bigg)\Bigg](x_0)d\xi_{m}\sigma(\xi')dx'\nonumber\\
&=&\frac{1}{2}\int_{|\xi'|=1}\int^{+\infty}_{-\infty}{\rm trace}
\Bigg[\partial^{2}_{\xi_{m}}\pi^{+}_{\xi_{m}}\sigma_{-1}(D^{-1})
      \times\partial_{x_{m}}\sigma_{-2n-1}
      \bigg(D^{-2n-1}\bigg)\Bigg](x_0)d\xi_{m}\sigma(\xi')dx'\nonumber\\
      &&-\frac{1}{2}f\cdot\partial_{x_{m}}(f^{-1})
      \int_{|\xi'|=1}\int^{+\infty}_{-\infty}{\rm trace}
\Bigg[\partial_{\xi_{m}}\pi^{+}_{\xi_{m}}\sigma_{-1}(D^{-1})
      \times\partial_{\xi_{m}}\sigma_{-2n-1}
\bigg(D^{-2n-1}\bigg)\Bigg](x_0)\nonumber\\
      &&\times d\xi_{m}\sigma(\xi')dx'.
\end{eqnarray}
By (2.2.29) in \cite{Wa3}, we have
\begin{eqnarray}\label{fum4.23}
\partial_{\xi_{m}}^2\pi^+_{\xi_{m}}\sigma_{-1}(D^{-1})(x_0)|_{|\xi'|=1}
=\frac{c(\xi')+ic(dx_m)}{(\xi_{m}-i)^3}.
\end{eqnarray}
By Lemma (\ref{lem3}), direct computations show that
\begin{eqnarray}\label{fum4.24}
\partial_{x_m}\sigma_{-1-2n}(D^{-1-2n})(x_0)|_{|\xi'|=1}=\frac{\sqrt{-1}\partial_{x_m}[c(\xi')](x_0)}
{(1+\xi_{m}^2)^{n+1}}-\frac{\sqrt{-1}(n+1)h'(0)c(\xi)}{(1+\xi^2_m)^{n+2}}.
\end{eqnarray}
According to the above three formulas and the Cauchy integral formula, we have
\begin{eqnarray}\label{fum4.25}
&&\frac{1}{2}\int_{|\xi'|=1}\int^{+\infty}_{-\infty}{\rm trace}
\Bigg[\partial^{2}_{\xi_{m}}\pi^{+}_{\xi_{m}}\sigma_{-1}(D^{-1})
      \times\partial_{x_{m}}\sigma_{-2n-1}\bigg(D^{-2n-1}\bigg)\Bigg](x_0)d\xi_{m}\sigma(\xi')dx'\nonumber\\
&=&\frac{1}{2}\int_{|\xi'|=1}\int^{+\infty}_{-\infty} {\rm trace}\left\{\frac{c(\xi')+ic(dx_m)}{(\xi_{m}-i)^3}\right.\left.\times
\left[\frac{i\partial_{x_m}[c(\xi')](x_0)}
{(1+\xi_{m}^2)^{n+1}}-\frac{i(n+1)h'(0)c(\xi)}{(1+\xi^2_m)^{n+2}}\right]\right\}
(x_0)d\xi_{m}\sigma(\xi')dx'\nonumber\\
&=&\frac{1}{2}\int_{|\xi'|=1}\int^{+\infty}_{-\infty}2^{n+1}h'(0)\times\frac{-i\xi^2_m-2(n+1)\xi_{m}+2ni}
{(\xi_{m}+i)^{n+2}(\xi_{m}-i)^{n+5}}d\xi_{m}\sigma(\xi')dx'\nonumber\\
&=&\frac{\pi ih'(0)2^{n+1}{\rm Vol}(S_{2n+2})dx'}{(n+4)!}G_1.
\end{eqnarray}
where $$G_1:=\left[\frac{-i\xi^2_m-2(n+1)\xi_{m}+2ni}{(\xi_{m}+i)^{n+2}}\right]^{(n+4)}|_{\xi_{m}=i}.$$
By (2.2.29) in \cite{Wa3}, we have
\begin{eqnarray}\label{fum4.26}
\partial_{\xi_{m}}\pi^+_{\xi_{m}}\sigma_{-1}
(D^{-1})(x_0)|_{|\xi'|=1}
=-\frac{c(\xi')+ic(dx_m)}{2(\xi_{m}-i)^2}.
\end{eqnarray}
Combining (\ref{fum4.17}) and (\ref{fum4.26}), we have
\begin{eqnarray}\label{fum4.27}
{\rm trace} \Big[\partial_{\xi_{m}}\pi^{+}_{\xi_{m}}\sigma_{-1}(D^{-1})
      \times\partial_{\xi_{m}}\big(\sigma_{-2n-1}(D^{-2n-1})\big)\Big](x_{0})|_{|\xi'|=1}
      =-\frac{2^{n+2}(n+1)i\xi_{m}(1+\xi_{m})}{(\xi_{m}-i)^{n+4}(\xi_{m}+i)^{n+2}},
\end{eqnarray}
then we obtain
\begin{eqnarray}\label{fum4.28}
&&-\frac{1}{2}f\cdot\partial_{x_{m}}(f^{-1})
      \int_{|\xi'|=1}\int^{+\infty}_{-\infty}{\rm trace}
\Bigg[\partial_{\xi_{m}}\pi^{+}_{\xi_{m}}\sigma_{-1}(D^{-1})
      \times\partial_{\xi_{m}}\sigma_{-2n-1}
\bigg(D^{-2n-1}\bigg)\Bigg](x_0)d\xi_{m}\sigma(\xi')dx'\nonumber\\
&=&\frac{f\cdot\partial_{x_{m}}(f^{-1})\pi i }{(n+3)!}\cdot {\rm Vol}(S_{2n+2})dx'\cdot G_2,
\end{eqnarray}
where $$G_2:=\bigg[\frac{2^{n+2}(n+1)i\xi_{m}(1+\xi_{m})}{(\xi_{m}+i)^{n+2}}\bigg]^{(n+3)}|_{\xi_{m}=i}.$$
Then
\begin{eqnarray*}
{\bf case~(3)}
&=&\Bigg[\frac{h'(0)2^{n+1}}{(n+4)!}G_1+\frac{f\cdot\partial_{x_{m}}(f^{-1})}{(n+3)!}G_2\Bigg]\cdot \pi i{\rm Vol}(S_{2n+2})dx'.
\end{eqnarray*}
\\
\noindent  {\bf case (4)}~$r=-1,l=-2n-2,|\alpha|=j=k=0$.\\

\noindent By (\ref{fum4.2}), we have
 \begin{eqnarray}\label{fum4.29}
&&{\rm case~ (4)}\nonumber\\
&=&-i\int_{|\xi'|=1}\int^{+\infty}_{-\infty}{\rm trace} \Big[\pi^{+}_{\xi_{m}}\sigma_{-1}(fD^{-1})
      \times\partial_{\xi_{m}}\sigma_{-2n-2}\bigg((f^{-1}D^{-1})\cdot (fD^{-1}\cdot f^{-1}D^{-1})^{n}\bigg)\Big](x_0)\nonumber\\
&&\times d\xi_{m}\sigma(\xi')dx'\nonumber\\
&=&i\int_{|\xi'|=1}\int^{+\infty}_{-\infty}{\rm trace} \Bigg[\partial_{\xi_{m}}\pi^{+}_{\xi_{m}}\sigma_{-1}(fD^{-1})
      \times\sigma_{-2n-2}\bigg((f^{-1}D^{-1})\cdot (fD^{-1}\cdot f^{-1}D^{-1})^{n}\bigg)\Bigg](x_0)\nonumber\\
&&\times d\xi_{m}\sigma(\xi')dx'\nonumber\\
&=&i\int_{|\xi'|=1}\int^{+\infty}_{-\infty}{\rm trace} \Bigg\{\partial_{\xi_{m}}\pi^{+}_{\xi_{m}}\big(f\sigma_{-1}(D^{-1})\big)
      \times  \bigg[(f^{-1}
      \sigma_{-2n-2}(D^{-2n-1})
      +\sum^{m}_{j=1}\partial_{\xi_{j}}(|\xi|^{-2n-2})\nonumber\\
  &&    \times\Big(\sigma_{1}(D)\partial_{x_{j}}(f)\Big)
      \bigg]\Bigg\} (x_0)d\xi_{m}\sigma(\xi')dx'\nonumber\\
&=&i\int_{|\xi'|=1}\int^{+\infty}_{-\infty}{\rm trace} \Big[\partial_{\xi_{m}}\pi^{+}_{\xi_{m}}\Big(\sigma_{-1}(D^{-1})\Big)\times\sigma_{-2n-2}(D^{-2n-1})\Big] (x_0)d\xi_{m}\sigma(\xi')dx'\nonumber\\
&&+i\int_{|\xi'|=1}\int^{+\infty}_{-\infty}{\rm trace} \Bigg[\partial_{\xi_{m}}\pi^{+}_{\xi_{m}}\big(\sigma_{-1}(D^{-1})\big)\times \sum^{m}_{j=1}\partial_{\xi_{j}}(|\xi|^{-2n-2})\Big[\sigma_{1}(D)\partial_{x_{j}}(f)\Big]\Bigg] (x_0)
\nonumber\\
&&\times d\xi_{m}\sigma(\xi')dx'.
\end{eqnarray}
By (3.8) in \cite{WJ5}, we have:
\begin{eqnarray}\label{fum5.1}
\sigma_{-2n-1}(D^{-2n})
=n\sigma_2(D^2)^{(-n+1)}\sigma_{-3}(D^{-2})-i \sum_{k=0}^{n-2}\sum_{\mu=1}^{2n+2}
\partial_{\xi_{\mu}}\sigma_{2}^{-n+k+1}(D^2)
\partial_{x_{\mu}}\sigma_{2}^{-1}(D^2)\big(\sigma_2(D^2)\big)^{-k}.
\end{eqnarray}
By Lemma 2.2 in \cite{Wa3}, we have
\begin{eqnarray}\label{fum5.2}
\sum_{j=1}^{m}\partial_{\xi_j}(|\xi|^{-2n-2})\partial_{x_j}(c(\xi))(x_0)|_{|\xi'|=1}
=-2(n+1)\xi_{m}(1+\xi^2_m)^{-n-2}\partial_{x_m}[c(\xi')](x_0),
\end{eqnarray}
and
\begin{eqnarray}\label{fum5.3}
\left[-i \sum_{k=0}^{n-1}\sum_{\mu=1}^{m}
\partial_{\xi_{\mu}}\sigma_{2}^{-n+k}(D^2)
\partial_{x_{\mu}}\sigma_{2}^{-1}(D^2)(\sigma_2(D^2))^{-k}\right]ic(\xi)(x_0)|_{|\xi'|=1}
=\frac{c(\xi)h'(0)\xi_{m}(n+1)n}{(1+\xi_{m}^2)^{n+3}}.
\end{eqnarray}
By (3.26) in \cite{WJ5}, we have
\begin{eqnarray}
\sigma_{-3}(D^{-2})(x_0)|_{|\xi'|=1}
&=&\frac{i}{(1+\xi_{m}^2)^2}\Big(\frac{1}{2}h'(0)\sum_{k<m}\xi_k c(\widetilde{e_k})c(\widetilde{e_m})
-\frac{n-1}{2}h'(0)\xi_{m}\Big)-\frac{2ih'(0)\xi_{m}}{(1+\xi_{m}^2)^3}\nonumber\\
&=&\frac{i}{(1+\xi_{m}^2)^2}\Big(\frac{1}{2}h'(0)c(\xi')c(dx_m)
                                       -\frac{n-1}{2}h'(0)\xi_{m}\Big)-\frac{2ih'(0)\xi_{m}}{(1+\xi_{m}^2)^3}.
\end{eqnarray}
So by (\ref{fum5.1}), we have
\begin{eqnarray}\label{fum5.4}
&&\sigma_{-2n-2}(D^{-2n-1})(x_0)|_{|\xi'|=1}=\sigma_{-2n-2}(D^{-2n-2}\cdot D)\nonumber\\
&=&\left\{\sum_{|\alpha|=0}^{+\infty}(-i)^{|\alpha|}
\frac{1}{\alpha!}\partial^\alpha_\xi[\sigma(D^{-2n-2})]\partial^\alpha_x[\sigma(D)]\right\}_{-2n-2}\nonumber\\
&=&\sigma_{-2n-2}(D^{-2n-2})\sigma_0(D)+\sigma_{-2n-3}(D^{-2n-2})\sigma_1(D)+\sum_{|\alpha|=1}(-i)
\partial^\alpha_\xi[\sigma_{-2n-2}(D^{-2n-2})]\partial^\alpha_x[\sigma_1(D)]\nonumber\\
&=&|\xi|^{-2n-2}\sigma_0(D)+\sum_{j=1}^{2n+4}\partial_{\xi_j}(|\xi|^{-2n-2})\partial_{x_j}c(\xi)+
\bigg[(n+1)\sigma_2(D^2)^{(-n)}\sigma_{-3}(D^{-2})-i \sum_{k=0}^{n-1}\sum_{\mu=1}^{2n+4}\partial_{\xi_{\mu}}\nonumber\\
&&\times\sigma_{2}^{-n+k}(D^2)
\partial_{x_{\mu}}\sigma_{2}^{-1}(D^2)(\sigma_2(D^2))^{-k}\bigg]\sqrt{-1}c(\xi)\nonumber\\
&=&\frac{(-2n-3)h'(0)c(dx_m)}{4(1+\xi_{m}^2)^{n+1}}
-2(n+1)\xi_{m}(1+\xi^2_m)^{-n-2}\partial_{x_m}
[c(\xi')](x_0)+(n+1)i(1+\xi_{m}^2)^{-n}[c(\xi')+\xi_m\nonumber\\
&&\times c(dx_m)]
\times\Bigg[
\frac{-ih'(0)c(\xi')c(dx_m)}{2(1+\xi_{m}^2)^2}-\frac{(2n+3)h'(0)i\xi_{m}}{2(1+\xi_{m}^2)^2}
-\frac{2ih'(0)\xi_{m}}{(1+\xi_{m}^2)^3}\Bigg]+[c(\xi')+\xi_mc(dx_m)]h'(0)\xi_{m}\nonumber\\
&&\times(n^2+n)(1+\xi_{m}^2)^{-n-3}.
\end{eqnarray}
By (\ref{fum4.26}) and (\ref{fum5.4}), we have
\begin{eqnarray}\label{fum5.5}
&&{\rm trace} [\partial_{\xi_{m}}\pi^+_{\xi_{m}}\sigma_{-1}(D^{-1})\times
        \sigma_{-2n-2}(D^{-2n-1})](x_0)|_{|\xi'|=1}\nonumber\\
&=&\frac{2^{n+1}h'(0)}
        {4(\xi_{m}-i)^{n+4}(\xi_{m}+i)^{n+3}}\times
        \bigg[(2n+3)(2n+1)i\xi_{m}^3+(2\pi-2n-1)\xi_{m}^2+\big(8n^2+16n+7\big)i\xi_{m}\nonumber\\
        &&+(2n+3+2\pi)\bigg]^{(n+3)}|_{\xi_n=i}
\end{eqnarray}
Then by the Cauchy integral formula, we get
\begin{eqnarray}\label{fum5.6}
&&i\int_{|\xi'|=1}\int^{+\infty}_{-\infty}{\rm trace} \Big[\partial_{\xi_{m}}\pi^{+}_{\xi_{m}}\big(\sigma_{-1}(D^{-1})\big)\times\sigma_{-2n-2}(D^{-2n-1})\Big] (x_0)d\xi_{m}\sigma(\xi')dx'\nonumber\\
&=&\frac{-\pi 2^{n}h'(0){\rm Vol }(S_{2n+2})dx'}{(n+3)!}G_3.
\end{eqnarray}
where $$G_3:=
\Bigg \{\frac{1} {(\xi_{m}+i)^{n+3}}
 \bigg[(2n+3)(2n+1)i\xi_{m}^3+(2\pi-2n-1)\xi_{m}^2+\big(8n^2+16n+7\big)i\xi_{m}+(2n+3+2\pi)\bigg]\Bigg\}^{(n+3)} |_{\xi_{m}=i}.$$
And we have
\begin{eqnarray}\label{fum5.7}
\sum_{j=1}^{m}\partial_{\xi_j}(|\xi|^{-2n-2})\bigg(\sigma_{1}(D)\partial_{x_j}(f)\bigg)(x_0)|_{|\xi'|=1}
=-\sum_{j=1}^{m}\bigg(\xi_j\partial_{x_j}(f)\bigg)\cdot(2n+2)i\cdot(\xi^{2}_m+1)^{-n-2}\cdot c(\xi).
\end{eqnarray}
We note that $i<m,~\int_{|\xi'|=1}\xi_i\sigma(\xi')=0$,
so $\sum\limits_{j}\xi_j\partial_{x_j}(f){\rm trace}[{\rm id}]$ have no contribution for computing {\bf case (4)}.
Then we obtain
\begin{eqnarray}\label{fum5.8}
&&i\int_{|\xi'|=1}\int^{+\infty}_{-\infty}{\rm trace} \Big[\partial_{\xi_{m}}\pi^{+}_{\xi_{m}}\big(\sigma_{-1}(D^{-1})\big)\times \sum^{m}_{j=1}\partial_{\xi_{j}}(|\xi|^{-2n-2})\big(\sigma_{1}(D)\partial_{x_{j}}(f)\big)\Big] (x_0)d\xi_{m}\sigma(\xi')dx'\nonumber\\
&=&i\int_{|\xi'|=1}\int^{+\infty}_{-\infty}
\bigg\{-
\frac{c(\xi')+ic(dx_m)}{2(\xi_{m}-i)^2}
\cdot\Big[-\Big(\xi_m\partial_{x_m}(f)\Big)\cdot(2n+2)i\cdot(\xi^{2}_m+1)^{-n-2}\cdot c(\xi)\Big]\bigg\}(x_0)\nonumber\\
&&\times d\xi_{m}\sigma(\xi')dx'\nonumber\\
&=&\frac{-2^{n+4}\partial_{x_{m}}(f)i\pi(n+1)}{(n+3)!}{\rm Vol}(S_{2n+2})G_4 dx'.
\end{eqnarray}
where$$G_4:=\bigg[\frac{(1+\xi_m)\xi_m}{2(\xi_m+i)^{n+2}}\bigg]^{(n+3)}|_{\xi_{m}=i}.$$
Then
\begin{eqnarray*}
{\bf case~(4)}=\bigg[\frac{-\pi 2^{n}h'(0)}{(n+3)!}
G_3-\frac{2^{n+4}\partial_{x_{m}}(f)i\pi(n+1)}{(n+3)!}G_4\bigg]{\rm Vol}(S_{2n+2})dx'.
\end{eqnarray*}

\noindent  {\bf case (5)}~$r=-2,~l=-2n-1,~k=j=|\alpha|=0.$\\

\noindent By (\ref{fum4.2}), we get
\begin{eqnarray}\label{fum5.9}
{\rm case~ (5)}&=&-i\int_{|\xi'|=1}\int^{+\infty}_{-\infty}{\rm trace} \bigg[\pi^+_{\xi_m}\sigma_{-2}(fD^{-1})\times
     \partial_{\xi_m}\sigma_{-2n-1}\Big((f^{-1}D^{-1})\cdot (fD^{-1}\cdot f^{-1}D^{-1})^{n}\Big)\bigg](x_0)\nonumber\\
     &&\times d\xi_m\sigma(\xi')dx'\nonumber\\
&=&-i\int_{|\xi'|=1}\int^{+\infty}_{-\infty}{\rm trace} \bigg[\pi^+_{\xi_m}\sigma_{-2}(D^{-1})\times
     \partial_{\xi_m}\sigma_{-2n-1}(D^{-2n-1})\bigg](x_0)d\xi_m\sigma(\xi')dx'.
\end{eqnarray}
\noindent By (2.2.34)-(2.2.37) in \cite{Wa3}, we have
\begin{eqnarray}\label{fum5.10}
\pi^+_{\xi_n}\sigma_{-2}(D^{-1})(x_0)|_{|\xi'|=1}
=J_1-J_2,
\end{eqnarray}
\noindent where\\
\begin{eqnarray}\label{fum5.11}
J_1=-\frac{H_1}{4(\xi_n-i)}-\frac{H_2}{4(\xi_n-i)^2},
\end{eqnarray}
\noindent and
\begin{eqnarray}\label{fum5.12}
H_1=ic(\xi')\sigma_0(D)c(\xi')+ic(dx_n)
\sigma_0(D)c(dx_n)+ic(\xi')c(dx_n)
\partial_{x_n}[c(\xi')];
\end{eqnarray}
\begin{eqnarray}\label{fum5.13}
H_2=[c(\xi')+ic(dx_n)]\sigma_0(D)[c(\xi')+ic(dx_n)]+c(\xi')c(dx_n)\partial_{x_n}c(\xi')-i\partial_{x_n}[c(\xi')];
\end{eqnarray}
\begin{eqnarray}\label{fum6.1}
J_2
=\frac{h'(0)}{2}\left[\frac{c(dx_n)}{4i(\xi_n-i)}+\frac{c(dx_n)-ic(\xi')}{8(\xi_n-i)^2}
+\frac{3\xi_n-7i}{8(\xi_n-i)^3}[ic(\xi')-c(dx_n)]\right].
\end{eqnarray}
Similar to (2.2.38) in \cite{Wa3}, we have
\begin{eqnarray}\label{fum5.14}
\partial_{\xi_m}\sigma_{-2n-1}(D^{-2n-1})(x_0)|_{|\xi'|=1}
=\sqrt{-1}\left[\frac{c(dx_m)}{(1+\xi_m^2)^{n+1}}-(n+1)\times\frac{2\xi_mc(\xi')+2\xi_m^2c(dx_m)}
{(1+\xi_m^2)^{n+2}}\right].
\end{eqnarray}
\noindent By (\ref{fum5.9}), (\ref{fum6.1}) and (\ref{fum5.14}), we have
\begin{eqnarray}\label{fum5.15}
&&{\rm tr }[J_2\times\partial_{\xi_m}\sigma_{-1-2n}(D^{-1-2n})(x_0)]|_{|\xi'|=1}\nonumber\\
&=&\frac{\sqrt{-1}}{2}h'(0){\rm trace}
\Bigg\{\bigg[\Big(\frac{1}{4i(\xi_m-i)}+\frac{1}{8(\xi_m-i)^2}-\frac{3\xi_m-7i}{8(\xi_m-i)^3}\Big)c(dx_m)
+\left(\frac{-1}{8(\xi_m-i)^2}+\frac{3\xi_m-7i}{8(\xi_m-i)^3}\right)\nonumber\\
&&\times ic(\xi')\bigg]\times\bigg[\left(\frac{1}{(1+\xi_m^2)^{1+n}}-\frac{2(n+1)\xi_m^2}{(1+\xi_m^2)^{n+2}}\right)c(dx_m)-\frac{2(n+1)\xi_m}
{(1+\xi_m^2)^{n+2}}c(\xi')\bigg]\Bigg\}\nonumber\\
&=&h'(0)2^{n-1}\times \frac{(2n+1)\xi^3_m-2i(2n+1)\xi^2_m-(6n+5)\xi_m+4i}{(\xi_m-i)^2(1+\xi^2_m)^{n+2}}.
\end{eqnarray}

\noindent By (2.2.40) in \cite{Wa3}, we have
\begin{eqnarray}\label{fum5.16}
J_1&=&\frac{-1}{4(\xi_m-i)^2}[(2+i\xi_m)c(\xi')\sigma_0(D)c(\xi')+i\xi_mc(dx_m)\sigma_0(D)c(dx_m)+(2+i\xi_m)c(\xi')c(dx_m)\nonumber\\
&&\times\partial_{x_m}c(\xi')+ic(dx_m)\sigma_0(D)c(\xi')
+ic(\xi')\sigma_0(D)c(dx_m)-i\partial_{x_m}c(\xi')].
\end{eqnarray}
Similar to Lemma 2.4 in \cite{Wa3}, we have
\begin{eqnarray}\label{fum8.1}
\sigma_0(D)(x_0)=c_0c(dx_m),~~{\rm where}~ c_0=\frac{-m}{4}h'(0).
\end{eqnarray}
\noindent By the relation of the Clifford action and ${\rm trace}{AB}={\rm trace }{BA}$, then we have the equalities:
$${\rm trace}[c(\xi')\sigma_0(D)c(\xi')c(dx_m)]=-c_02^{n+2};~~{\rm trace}[c(dx_m)\sigma_0(D)c(dx_m)^2]=c_02^{n+2};$$
$${\rm trace}[c(\xi')c(dx_m)\partial_{x_m}c(\xi')c(dx_m)](x_0)|_{|\xi'|=1}=-2^{n+1}h'(0);~
{\rm trace}[c(dx_m)\sigma_0(D)c(\xi')^2]=c_02^{n+2}.$$

\noindent By (4.47) and (4.48),
 considering for $i<m$, $\int_{|\xi'|=1}\{{\rm odd ~number~ product~ of~}\xi_i\}\sigma(\xi')=0$, then
\begin{eqnarray}\label{fum8.2}
&&{\rm tr }[J_1\times\partial_{\xi_m}\sigma_{-1-2n}(D^{-1-2n})(x_0)]|_{|\xi'|=1}\nonumber\\
&=&\frac{2^{n+1}ih'(0)}{4(\xi_m-i)^2(1+\xi_m^2)^{n+2}}\cdot\bigg\{
(m-1)[(2n+1)\xi_m^2-2i(n+1)\xi_m-1]+[-(1+2n)i\xi_m^3-2(1+2n)\xi_m^2\nonumber\\
&&+(2n+3)i\xi_m+2]\bigg\}.
\end{eqnarray}
By combining (\ref{fum5.15}), (\ref{fum8.2}) and the Cauchy integral formula, we have
\begin{eqnarray}\label{fum7.2}
&&{\rm case~(5)}\nonumber\\
&=&-i\int_{|\xi'|=1}\int^{+\infty}_{-\infty}{\rm trace} [(J_1-J_2)\times
     \partial_{\xi_m}\sigma_{-1-2n}(D^{-1-2n})](x_0)d\xi_m\sigma(\xi')dx'\nonumber\\
 &=&2^{n+1}h'(0)\int_{|\xi'|=1}\int^{+\infty}_{-\infty}\frac{[4n^2+8n+3]\xi_m^2
 -[4n^2+14n+8]i\xi_m-(m+1)}{4(\xi_m-i)^{n+4}(\xi_m+i)^{n+2}}
     d\xi_m\sigma(\xi')dx'\nonumber\\
&=&\frac{2^{n+2}h'(0){\rm Vol}(S_{2n+2})\pi i dx'}{(n+3)!}G_5,
\end{eqnarray}
where $$G_5:=\left[
    \frac{[4n^2+8n+3]\xi_m^2-[4n^2+14n+8]i\xi_m-(m+1)}{4(\xi_m+i)^{n+2}}\right]
    ^{(n+3)}|_{\xi_m=i}. $$
Since $\Psi$ is the sum of the {\bf case~(1)-case~(5)}, so
\begin{eqnarray}\label{fum10.1}
\Psi&=&\Bigg\{\frac{(-1)^{n}h'(0)\pi }{3\times2^{n+6}(3+n)!}Y_0
+\frac{(-1)^{n}(n+1)f^{-1}\partial_{x_{n}}(f)\pi}{2^{n+2}}Y_1
+\bigg[\frac{(i-1)f\cdot\partial_{x_{n}}(f^{-1})}{2^{n+3}}\nonumber\\
&&-\frac{\partial_{x_{n}}(f)(1+i)}{2^{n+2}} \bigg]Y_2\Bigg\}{\rm Vol}(S_{2n+2})dx',
\end{eqnarray}
where
\begin{eqnarray}\label{fum11.1}
Y_0&=&-24(1+i)\bigg[(4n^2+8n+3)C_{-n-3}^{n}+(6n^2+13n+5-\pi)C_{-n-3}^{n+1}+(n^{2}+3n-\pi+1)C_{-n-3}^{n+2}   \nonumber\\
&&-2n(1+i)(1+n)C_{-n-3}^{n+3}+24C_{-n-2}^{n+2}-24C_{-n-2}^{n+3}-6C_{-n-2}^{n+4}\bigg]\times (3+n)!+A^{3+n}_{-1-n};
\nonumber\\
Y_1&=&(4n+4+6i)C_{-n-3}^{n}+(9i+8n+9)C_{-n-3}^{n+1}+(i+n+1)\bigg[5C_{-n-3}^{n+2}+C_{-n-3}^{n+3}\bigg]
\nonumber\\
Y_2&=&(n+1)(-1)^{n}\bigg[2(1+i)C_{-n-2}^{n+1}+(3+i)C_{-n-2}^{n+2}+C_{-n-2}^{n+3} \bigg].
\end{eqnarray}
Combining (\ref{fum3.10}) and (\ref{fum10.1}), we obtain Theorem 1.2.

\section*{Acknowledgements}
This work was supported by NSFC No.12301063 and NSFC No.11771070, DUFE202159 and Basic research Project of the Education Department of Liaoning Province (Grant No. LJKQZ20222442). The authors thank the referee for his (or her) careful reading and helpful comments.

\end{document}